\documentclass[11 pt]{amsart}
\usepackage{amssymb}
\usepackage{amsmath}
\usepackage{amsthm}
\usepackage{amsfonts}
\usepackage{standalone}

\usepackage[margin=1 in]{geometry}
\usepackage{multicol}
\usepackage{tikz}
\usetikzlibrary{patterns}
\usepackage{lmodern,tikz}
\usetikzlibrary{positioning}
\usepackage{pgfplots,caption}
\pgfplotsset{compat=1.13}
\usepackage[utf8]{inputenc}
\usepackage{graphicx}
\usepackage{verbatim}
\usepackage{enumerate}
\usepackage{mathtools}
\usepackage{xcolor}
\usetikzlibrary{calc}
\usepackage{float}

\usepackage{mathpazo}
\usepackage{euler}
\usepackage{microtype}

\DeclarePairedDelimiter\abs{\lvert}{\rvert}%
\DeclarePairedDelimiter\norm{\lVert}{\rVert}%

\makeatletter
\let\oldabs\abs
\def\abs{\@ifstar{\oldabs}{\oldabs*}}
\let\oldnorm\norm
\def\norm{\@ifstar{\oldnorm}{\oldnorm*}}
\makeatother

\usepackage{amsmath,stackengine,scalerel}
\stackMath
\def\dhat#1{\ThisStyle{\setbox0=\hbox{$\SavedStyle#1$}%
  \stackengine{0pt}{\SavedStyle#1}{\SavedStyle\hspace{.4\ht2}%
  \hat{\vphantom{#1}}\kern\dimexpr3.2\LMpt+2.0pt\relax\hat{\vphantom{#1}}}{O}{c}{F}{T}{L}}%
}

\tikzset{
    triple/.style args={[#1] in [#2] in [#3]}{
        #1,preaction={preaction={draw,#3},draw,#2}
    }
}

\DeclareMathAlphabet{\mathcal}{OMS}{cmsy}{m}{n}

\newcommand{\R}{\mathcal{R}}
\newcommand{\B}{\mathcal{B}}
\newcommand{\M}{\mathcal{M}}
\newcommand{\G}{G}
\newcommand{\T}{\mathcal{T}}

\DeclarePairedDelimiter\parentheses{\lparen}{\rparen}
\newcommand{\type}[1]{\operatorname{type} \parentheses{#1}}
\newcommand{\GG}{\G_{\text{red}}}

\newcommand{\newword}[1]{\textbf{\emph{#1}}}

\usepackage{xpatch}
\makeatletter
\xpatchcmd{\@thm}{\thm@headpunct{.}}{\thm@headpunct{}}{}{}
\makeatother
\newtheorem{thm}{Theorem}[section]
\newtheorem{theorem}[thm]{Theorem}

\newtheorem{proposition}[thm]{Proposition}
\newtheorem{lemma}[thm]{Lemma}

\newtheorem{conjecture}[thm]{Conjecture}

\newtheorem{corollary}[thm]{Corollary}

\newtheorem{definition}[thm]{Definition}

\newtheorem{remark}[thm]{Remark}

\newtheoremstyle{case}{}{}{}{}{}{:}{ }{}
\theoremstyle{case}

\tikzset{
    single/.style args={}{
        draw,line width=.5mm,black
    }
}
\tikzset{
    double/.style args={}{
        preaction={draw,line width=1.25mm,black},draw,line width=.25mm,white
    }
}
\tikzset{
    triple/.style args={}{
        preaction={preaction={draw,line width=2mm,black},draw,line width=1mm,white},draw,line 
        width=.5mm,black
    }
}  
\tikzset{
    quadruple/.style args={}{
        preaction={preaction={preaction={draw,line width=2.75mm,black},draw,line width=1.75mm,white},draw,line width=1.25mm,black},draw,line width=.25mm,white
    }
}

\title{Triconed Graphs, weighted forests, and $h$-vectors of matroid complexes }
\author[J. David]{Jacob David}
\author[P. Lai]{Pierce Lai}
\author[S. Oh]{SuHo Oh}
\author[C. Wu]{Christopher Wu}

\address{Phillips Exeter Academy}
\email{jdavid@exeter.edu}

\address{Massachusetts Institute of Technology, Department of Electrical Engineering and Computer Science}
\email{pwlai@mit.edu}

\address{Texas State University, Department of Mathematics}
\email{s\_o79@txstate.edu}

\address{Westlake High School}
\email{cw86459@gmail.com}

\date{\today}

\begin{document}


\begin{abstract}
A well-known conjecture of Stanley is that the $h$-vector of a matroid is a pure ${\mathcal O}$-sequence. There have been numerous papers with partial progress on this conjecture, but it is still wide open. In particular, for graphic matroids coming from taking the spanning trees of a graph as bases, the conjecture is mostly unsolved. In graph theory, a set of vertices is called dominating if every other vertex is adjacent to some vertex inside the chosen set. Kook proved Stanley's conjecture for coned graphs, which is the class of graphs that are dominated by a single vertex. Cranford et al extended that result to biconed graphs, which is the class of graphs dominated by a single edge. In this paper we extend that result to triconed graphs, the class of graphs dominated by a path of length $2$.


\end{abstract}

\maketitle

\section{Introduction}
Given a graph $G$, its \newword{graphic matroid} can be thought of as the collection of spanning trees of $G$. A matroid is a structure which generalizes the notion of independence that rises in linear algebra, graph theory and other areas. Given a matroid, we can naturally construct a simplicial complex and study its so called \newword{h-vector} which provides the topological information of the complex. The following conjecture by Stanley has motivated numerous research on h-vectors of a matroid:
\begin{conjecture}[\cite{Stanley1977}]
\label{con:stanley}
The h-vector of  a matroid is a pure O-sequence.
\end{conjecture}

Here an O-sequence is a sequence that can be obtained from counting monomials of each degree in an \newword{order ideal}: a finite collection $X$ of monomials such that, whenever $M \in X$ and $N$ divides $M$, we have $N \in X$ as well. If all maximal monomials of $X$ have the same degree, we say it is \newword{pure} and the degree sequence coming from it a pure O-sequence.
The above conjecture has been open for over four decades, and there are some specific classes of matroids for which the above conjecture has been established to be true. In particular, these include cographic matroids by Merino in \cite{Merino2001}, lattice-path matroids by Schweig in \cite{Schweig2010}, cotransversal matroids by Oh in \cite{Oh2013}, paving matroids by Merino, Noble, Ramirez-Ibanez, and Villarroel-Flores  \cite{MNRV2012}, internally perfect matroids by Dall in \cite{Dal}, rank $3$ matroids by H\'a, Stokes, and Zanello in \cite{HaStokesZanello2013}, rank $3$ and corank $2$ matroids by DeLoera, Kemper, and Klee in \cite{DelKemKle}, rank $4$ matroids by Klee and Samper in \cite{KleeSamper2015}, and rank $d$ matroids with $h_d \leq 5$  by Constantinescu, Kahle, and Varbaro in \cite{ConKahVar}.

In this paper the focus will solely be on the class of graphic matroids, for which the conjecture is still wide open. In \cite{Kook2012} Kook proved the conjecture for the class of \newword{coned graphs}, which is the class of graphs dominated by a vertex: there is a conal vertex adjacent to every other vertex. The result was extended in \cite{biconed} to the class of \newword{biconed graphs}, which is the class of graphs having a dominating edge: every other vertices of the graph is adjacent to at least one vertex of that edge. In this paper, we will further extend this result to \newword{triconed graphs}, which is the class of graphs being dominated by a path of length two: every other vertices of the graph is adjacent to at least one vertex along this path.

As was the case for the proof of coned graphs in \cite{Kook2012} and biconed graphs in \cite{biconed}, the main idea will be to convert spanning trees to forests with certain edges being \newword{empowered}, which then directly translates to monomials where variables correspond to edges of the graph. In section 2 we recall some basic notion and definitions from matroid and graph theory. In section 3 we describe our main objects in consideration and establish bijections between them: spanning trees of a triconed graph $G$, trirooted forests of $G \setminus B_0$ and $3$-(edge)weighted forests of $G \setminus B_0$. In section $4$, we show that our construction results in a pure multicomplex and hence proves Stanley's conjecture for graphic matroids coming from triconed graphs. In section 5 we provde a detailed example to demonstrate our method and discuss further directions of research.

\section{Preliminary}

Throughout the paper, we will use the following notation. When $A$ is a set and $a$ is a single element, instead of writing $A \setminus \{a\}$ we will instead write $A \setminus a$. Similarly, instead of writing $A \cup \{a\}$ we will instead use $A \cup a$.

\subsection{Graphs and Matroids} 
We first review some basic notions and tools of matroid theory. We will assume basic familiarity with matroid theory and recommend \cite{BrylawskiOxley1992} for more details. Let us start with a graph $G$. The edge set $E$ is called the \newword{ground set} of the matroid $\M_G$ we are about to construct from $G$. Take the collection of subsets of $E$ that corresponds to spanning trees of $G$: such sets are called the \newword{bases} of the \newword{graphic matroid} of $\M_G$ and will be denoted by $\B(G)$ or just simply $\B$. All spanning trees have the same number of edges in them and that cardinality $r$ is called the \newword{rank} of the matroid.


The bases of a matroid satisfy the following \newword{exchange property}: for any pair of bases $A,B \in \B$, for all $a \in A \setminus B$ there exists some $b \in B \setminus A$ such that $A \setminus a \cup b \in \B$.

An element of the ground set $E$ that is not contained in any base of $\M$ is called a \newword{loop}. An element of the ground set $E$ that is contained in every base of $\M$ is called a \newword{coloop}. Two elements $a,b$ of the ground set $E$, such that we have the property that $I\cup a \in \B$ iff $I\cup b \in \B$ for all subsets $I \subset E$, are called \newword{parallel}.

\subsection{Activity and the h-vector} 
From any ordering on the ground set $E$, we get an induced lexicographic ordering on the bases. An element $i$ of a base $B$ is \newword{(internally) active} if $B \setminus i \cup j$ is not a base for any $j < i$. Otherwise, it is called \newword{(internally) passive}. There is also the notion of externally active and externally passive elements, but in this paper we will only stick to interally active and passive elements, and hence will skip the usage of the word internally throughout. Given a base $B$ of a matroid, we call its \newword{passivity} as the number of (internally) passive elements it contains.

We will use this activity property of h-vectors as the definition of h-vectors of matroids in this paper.

\begin{theorem}
Let $(h_0,\dots,h_r)$ be the $h$-vector of a matroid $\M$. For $0 \leq i \leq r$, the entry $h_i$ is the number of bases of $\M$ with $i$ passive elements.
\end{theorem}

\begin{figure}[h]

\begin{center}
            \begin{tikzpicture}[scale=.4]
        \node [draw,fill,circle,inner sep = 0pt, minimum size = .3cm] (0) at (-3,3) {};
        \node [draw,fill,circle,inner sep = 0pt, minimum size = .3cm] (1) at (3,3) {};
        \node [draw,fill,circle,inner sep = 0pt, minimum size = .3cm] (2) at (-3,-3) {};
        \node [draw,fill,circle,inner sep = 0pt, minimum size = .3cm] (3) at (3,-3) {};

        \draw [line width = .5mm, black](0) -- (1);
        \draw [line width = .5mm, black](0) -- (2);
        \draw [line width = .5mm, black](0) -- (3);  
        \draw [line width = .5mm, black](1) -- (3);  
        \draw [line width = .5mm, black](2) -- (3);  
        
        \def \shif{6cm}
        \def \shify{-6cm}

        \node [draw=white,fill=white,circle,inner sep = 0pt, minimum size = .5cm,label=center:{\large1}] at (0,3) {};
        \node [draw=white,fill=white,circle,inner sep = 0pt, minimum size = .5cm,label=center:{\large2}] at (-3,0) {};
        \node [draw=white,fill=white,circle,inner sep = 0pt, minimum size = .5cm,label=center:{\large3}] at (0,0) {};
        \node [draw=white,fill=white,circle,inner sep = 0pt, minimum size = .5cm,label=center:{\large4}] at (3,0) {};
        \node [draw=white,fill=white,circle,inner sep = 0pt, minimum size = .5cm,label=center:{\large5}] at (0,-3) {};
        
        
        \node [draw,fill,circle,inner sep = 0pt, minimum size = .2cm] (a0) at (5,5) {};
        \node [draw,fill,circle,inner sep = 0pt, minimum size = .2cm] (a1) at (9,5) {};
        \node [draw,fill,circle,inner sep = 0pt, minimum size = .2cm] (a2) at (5,1) {};
        \node [draw,fill,circle,inner sep = 0pt, minimum size = .2cm] (a3) at (9,1) {};
        
        \draw [line width = .5mm, black](a0) -- (a1);
        \draw [line width = .5mm, black](a0) -- (a2);
        \draw [line width = .5mm, black](a0) -- (a3);  
        
        
        \begin{scope}[xshift=\shif]

        \node [draw,fill,circle,inner sep = 0pt, minimum size = .2cm] (a0) at (5,5) {};
        \node [draw,fill,circle,inner sep = 0pt, minimum size = .2cm] (a1) at (9,5) {};
        \node [draw,fill,circle,inner sep = 0pt, minimum size = .2cm] (a2) at (5,1) {};
        \node [draw,fill,circle,inner sep = 0pt, minimum size = .2cm] (a3) at (9,1) {};
        
        \draw [line width = .5mm, black](a0) -- (a1);
        \draw [line width = .5mm, black](a0) -- (a2);
        \draw [line width = .5mm, red](a1) -- (a3);  
        
        \end{scope}

        
        \begin{scope}[xshift=2*\shif]
        
        \node [draw,fill,circle,inner sep = 0pt, minimum size = .2cm] (a0) at (5,5) {};
        \node [draw,fill,circle,inner sep = 0pt, minimum size = .2cm] (a1) at (9,5) {};
        \node [draw,fill,circle,inner sep = 0pt, minimum size = .2cm] (a2) at (5,1) {};
        \node [draw,fill,circle,inner sep = 0pt, minimum size = .2cm] (a3) at (9,1) {};
        
        \draw [line width = .5mm, black](a0) -- (a1);
        \draw [line width = .5mm, black](a0) -- (a2);
        \draw [line width = .5mm, red](a2) -- (a3); 

        \end{scope}
        
        
        \begin{scope}[xshift=3*\shif]
        
        \node [draw,fill,circle,inner sep = 0pt, minimum size = .2cm] (a0) at (5,5) {};
        \node [draw,fill,circle,inner sep = 0pt, minimum size = .2cm] (a1) at (9,5) {};
        \node [draw,fill,circle,inner sep = 0pt, minimum size = .2cm] (a2) at (5,1) {};
        \node [draw,fill,circle,inner sep = 0pt, minimum size = .2cm] (a3) at (9,1) {};
        
        \draw [line width = .5mm, black](a0) -- (a1);
        \draw [line width = .5mm, red](a0) -- (a3);  
        \draw [line width = .5mm, red](a2) -- (a3); 

        \end{scope}

        
        \begin{scope}[yshift=\shify]
        
        \node [draw,fill,circle,inner sep = 0pt, minimum size = .2cm] (a0) at (5,5) {};
        \node [draw,fill,circle,inner sep = 0pt, minimum size = .2cm] (a1) at (9,5) {};
        \node [draw,fill,circle,inner sep = 0pt, minimum size = .2cm] (a2) at (5,1) {};
        \node [draw,fill,circle,inner sep = 0pt, minimum size = .2cm] (a3) at (9,1) {};
        
        \draw [line width = .5mm, black](a0) -- (a1);
        \draw [line width = .5mm, red](a1) -- (a3);  
        \draw [line width = .5mm, red](a2) -- (a3); 
        
        \end{scope}

        
        \begin{scope}[yshift=\shify, xshift=\shif]
        
        \node [draw,fill,circle,inner sep = 0pt, minimum size = .2cm] (a0) at (5,5) {};
        \node [draw,fill,circle,inner sep = 0pt, minimum size = .2cm] (a1) at (9,5) {};
        \node [draw,fill,circle,inner sep = 0pt, minimum size = .2cm] (a2) at (5,1) {};
        \node [draw,fill,circle,inner sep = 0pt, minimum size = .2cm] (a3) at (9,1) {};
        
        \draw [line width = .5mm, black](a0) -- (a2);
        \draw [line width = .5mm, red](a0) -- (a3);  
        \draw [line width = .5mm, red](a1) -- (a3);  
        
        \end{scope}
        
        
        \begin{scope}[yshift=\shify, xshift=2*\shif]
        
        \node [draw,fill,circle,inner sep = 0pt, minimum size = .2cm] (a0) at (5,5) {};
        \node [draw,fill,circle,inner sep = 0pt, minimum size = .2cm] (a1) at (9,5) {};
        \node [draw,fill,circle,inner sep = 0pt, minimum size = .2cm] (a2) at (5,1) {};
        \node [draw,fill,circle,inner sep = 0pt, minimum size = .2cm] (a3) at (9,1) {};
        
        \draw [line width = .5mm, red](a0) -- (a2);
        \draw [line width = .5mm, red](a1) -- (a3);  
        \draw [line width = .5mm, red](a2) -- (a3); 
        
        \end{scope}
        
        
        \begin{scope}[yshift=\shify, xshift=3*\shif]
        
        \node [draw,fill,circle,inner sep = 0pt, minimum size = .2cm] (a0) at (5,5) {};
        \node [draw,fill,circle,inner sep = 0pt, minimum size = .2cm] (a1) at (9,5) {};
        \node [draw,fill,circle,inner sep = 0pt, minimum size = .2cm] (a2) at (5,1) {};
        \node [draw,fill,circle,inner sep = 0pt, minimum size = .2cm] (a3) at (9,1) {};
        
        \draw [line width = .5mm, red](a0) -- (a3);  
        \draw [line width = .5mm, red](a1) -- (a3);  
        \draw [line width = .5mm, red](a2) -- (a3); 
        
        \end{scope}
        
    \end{tikzpicture}
\captionsetup{width=0.85\linewidth}
  \caption{An example of a graphic matroid. Passive elements are denoted by red. The h-vector is thus given by $\{1,2,2,3\}$.}
  \label{fig:exgraph1}
\end{center}

\vfill\null

\end{figure}
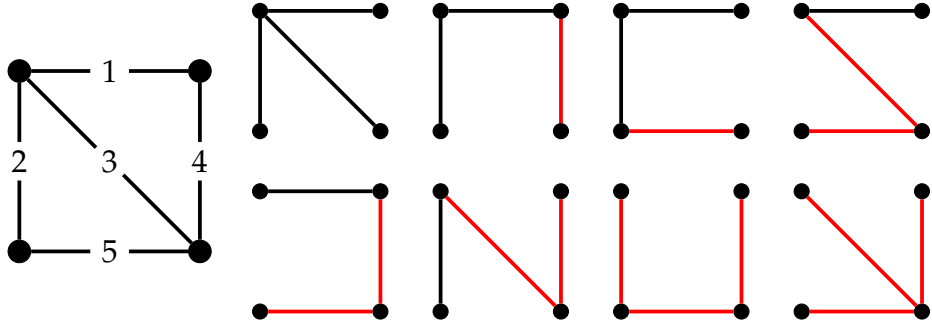

For example, take a look at Figure~\ref{fig:exgraph1}. The graph $G$ in consideration has $5$ edges, labeled using the set $[5] = \{1,\ldots,5\}$. There are $8$ spanning trees of $G$ in total. If we take the edge set of each spanning tree, we get $\B(G) = \{123,124,125,135,145,234,245,345\}$ where $abc$ is an abbreviation for $\{a,b,c\}$. They form the set of bases of the graphic matroid coming from $G$. Take a look at $125$. Here the edge $5$ is passive, since we can replace it with a smaller edge (say $3$) to get another spanning tree. Note that $5$ was colored red: in this figure, for each spanning tree, the passive edges are colored red. Hence the $h$-vector of the graphic matroid will be $(1,2,3,2)$.

\subsection{Order ideal and pure O-sequences}
An \newword{order ideal} is a finite collection $X$ of monomials such that, whenever $M \in X$ and $N$ divides $M$, we have $N \in X$ as well. Given an order ideal its \newword{O-sequence} is a sequence that can be obtained from counting monomials of each degree. We say that the order ideal is \newword{pure} if all the maximal monomials are of the same degree, and also call the O-sequence pure if it comes from a pure order ideal.

For example, take a look at Figure~\ref{fig:order_ideal}. If we take any monomial among this collection, say $xy^2$, we can see that all divisors of it is again present in the collection. Hence this is an order ideal. Moreover, since the maximal monomials under division (which are usually not guaranteed to all have the same degree) here are $x^3$ and $xy^2$, which have the same degree, this is a pure order ideal. Now notice that the degree sequence of this pure order ideal is $(1,2,3,2)$. Hence we can see that the $h$-vector of the graphic matroid back in Figure~\ref{fig:exgraph1}, which was $(1,2,3,2)$, is indeed a pure $O$-sequence. Hence we have just verified that Stanley's conjecture is true for the graph in Figure~\ref{fig:exgraph1}. What we want to do in this paper is to similarly construct a pure order ideal for the class of triconed graphs, hence proving the conjecture for this class.

\begin{figure}
    \centering
    \begin{tikzpicture}
        \node at (-1.5,2.4) {$x^3$};
        \node at (0.5,2.4) {$xy^2$};
        \node at (-1,1.6) {$x^2$};
        \node at (0,1.5) {$xy$};
        \node at (1,1.6) {$y^2$};
        \node at (-0.5,0.8) {$x$};
        \node at (0.5,0.8) {$y$};
        \node at (0,0) {$1$};
    \end{tikzpicture}
    \caption{Example of a pure order ideal. The O-sequence of this is $\{1,2,3,2\}$.}
    \label{fig:order_ideal}
\end{figure}


\section{Triconed graphs and rooted forests} \label{sec:triconed}

In this section we will go over the core definitions and ideas needed for this paper on triconed graphs.  The rough sketch of the strategy we use for attacking Stanley's conjecture for graphs is the following: we will show an algorithm for converting spanning trees into monomials, that will turn out to form a pure order ideal. The variables we use will come from edges of the graph $G$ that are not contained in a certain canonical spanning tree $B_0$ which will be defined in the subsection below. So given a spanning tree $B$, we convert the information of edges in $B_0 \cap B$ to marks on certain vertices, then convert the information stored on the marks of vertices into empowerment of edges which corresponds to increasing the power of the corresponding variable.

\subsection{Triconed graphs}
We shall define triconed graphs as the following:

\begin{definition}\label{tricone}
Given a graph $G$, we call $G$ a \newword{triconed graph} if one can choose vertices $0$, $1$, and $2$, with the following conditions: 
\begin{itemize}
    \item $1$ and $2$ are adjacent to $0$,
    \item there are no parallel edges having $0$ or $1$ or $2$ as endpoints,
    \item all other vertices are adjacent to at least $1$ of those $3$ vertices. 
    \end{itemize}
Alternatively, we may say that a triconed graph is a graph which is dominated by a path of length $2$, and which has no parallel edges touching this path.
\end{definition}

We call the vertices $0,1,2$ to be \newword{special} and other remaining vertices to be \newword{normal}.

\begin{remark}
\label{rem:conedandbiconed}
Coned graphs from \cite{Kook2007} and biconed graphs from \cite{biconed} are special cases of triconed graphs as long as they contain a pair of adjacent edges. Also the restriction of possible parallel edges is in line with that from \cite{biconed}: a graph is a biconed graph if it is dominated by an edge and no parallel edges touch this edge.
\end{remark}

The goal of our paper is to prove Stanley's conjecture for triconed graphs just defined. We start by showing that loops and coloops can be ignored:

\begin{lemma}
\label{lem:nonewparallel}
Let $G$ be a triconed graph. If we delete a loop (an edge not contained in any spanning tree) or contract by a coloop (take an edge contained in every spanning tree and merge the two endpoints) we get a biconed graph or a triconed graph.
\end{lemma}
\begin{proof}
The statement regarding the loops is obvious. If we have an edge $ab$ as a coloop inside the graph, we cannot have both $0a,0b$ inside the graph. Since otherwise, we can come up with a spanning tree containing $0a$ and $0b$ but not $ab$ using Kruskal's algorithm \cite{10.2307/2033241}. The same occurs for pairs $1a,1b$ and $2a,2b$. Hence after contracting by the edge $ab$, we avoid creating new parallel edges having a special vertex as an endpoint. 
\end{proof}

From the above lemma we can conclude the following:

\begin{remark}
\label{rem:noloopcoloop}
Given a triconed graph, deleting a loop or contracting by a coloop does not change the $h$-vector of the matroid. Hence for sake of proving Stanley's conjecture on the class of triconed graphs, we may assume that our graph does not have any loops or coloops thanks to Lemma~\ref{lem:nonewparallel}. Hence from now on, throughout the entire paper, we are only going to deal with triconed graphs with no loops nor coloops.
\end{remark}

We will put an ordering on the vertices and then the edges. First for the vertices, label $0,1,2$ as above. Then we label the remaining vertices adjacent to $0$ starting from $3$. After that, we label the remaining vertices adjacent to $1$ then finish off with the remaining vertices adjacent to $2$. This gives a total ordering on the vertices of $G$. For example look at Figure~\ref{fig:exgraph2}. We have vertex $0$ that is adjacent to $1,2$ and all other vertices are adjacent to at least one of $0,1,2$. The remaining vertices adjacent to $0$ are labeled $3,4$, the remaining vertices adjacent to $1$ are labeled $5,6,7$ and finally the remaining vertices (which are adjacent to $2$) are labeled $8,9,A$ (in this paper we will use $A$ when used as a label of a vertex, to denote $10$). We will call the vertices adjacent to $0$ as vertices of \newword{height 1} of the graph and the vertices not adjacent to $0$ (also excluding $0$ as well) as vertices of \newword{height 2} of the graph.

Now on to the total ordering on the edges. The smallest edges will be the ones adjacent to $0$. Among that collection we order them lexicographically. The next smallest edges will be the ones connecting $1$ to vertices of height $2$. Within that collection of edges we again order them lexicographically. Then the next set of edges will be the ones connecting $2$ to vertices of height $2$ that are not adjacent to $1$. Again order within lexicographically. Finally comes the set of all the remaining edges, which we can order in any way. For example take a look at Figure~\ref{fig:exgraph2}. Then the ordering on the edges will be as the following (here the letter $A$ stands for $10$):
$$01,02,03,04,15,16,17,28,29,2A,25,27,34,35,36,38,47,49,56,78,9A.$$

Notice that in the above order, we are free to reorder the edges $25,27,\ldots,9A$ in any way we want. 
Using the total ordering on the edges we just went over, we define the \newword{canonical spanning tree} $B_0$ of $G$ as the lexicographically smallest spanning tree of $G$. For example, in Figure~\ref{fig:exgraph2}, the lexicographically smallest spanning tree consists of the edges $01,02,03,04,15,16,17,28,29,2A$ (the corresponding edges are colored red). Another way to think about this canonical spanning tree is that it is obtained by doing a breadth-first search starting from vertex $0$, then $1$ and $2$. Note that since $B_0$ is the lexicographically smallest spanning tree of $\G$, it has the property that no edges are passive.


We say that vertex $u$ is a child of vertex $v$ if $u$ is a child of $v$ in the canonical spanning tree $B_0$ rooted at vertex $0$. For an edge $vu$ of $B_0$ where $u$ is a child of $v$, we say that the child of the edge $vu$ is $u$ and the \newword{cone edge} of $u$ is $vu$. Notice that this gives us a bijection between non-root vertices of $B_0$ and edges of $B_0$: simply map a vertex to its cone edge. For example, take a look at the canonical spanning tree drawn in red inside Figure~\ref{fig:exgraph2}. The cone edge of vertex $5$ would be $15$, and the cone edge of the vertex $3$ would be $03$.

We classify a vertex as $\newword{type 1}$ if it is vertex $1$ or a child of vertex $1$, as $\newword{type 2}$ if it is vertex $2$ or a child of vertex $2$, and as $\newword{type 0}$ otherwise. We will look mainly at subgraphs resulting from the removal of the spanning tree $B_0$. Let $\GG$ be the graph obtained from $G$ by removing the edges of $B_0$ and the vertex $0$. For example take a look at the triconed graph $G$ in Figure~\ref{fig:exgraph2}. Removing $0$ and edges of $B_0$ from this graph gives us $\GG$ in Figure~\ref{fig:exgraphcone}. In $\GG$, vertices $1,5,6,7$ are of type $1$, vertices $2,8,9,A$ are of type $2$ and vertices $3,4$ are of type $0$.

\begin{figure}[h]
\begin{multicols}{2}

\begin{center}
    \tikzset{
        unmarked/.style args={}{
            draw,fill = black,circle,inner sep = 0pt, minimum size = .30cm,
        }
    }
    \begin{tikzpicture}[scale=1]
        \node [unmarked={},label=below:{\small5}] (M1) at (-2.5,0) {};
        \node [unmarked={},label=below:{\small6}] (M2) at (-1.5,0) {};
        \node [unmarked={},label=below:{\small7}] (M3) at (-0.5,0) {};
        \node [unmarked={},label=below:{\small8}] (M4) at (0.5,0) {};
        \node [unmarked={},label=below:{\small9}] (M5) at (1.5,0) {};
        \node [unmarked={},label=below:{\small10}] (M6) at (2.5,0) {};
        \node [unmarked={},label=above:{\small3}] (N1) at (0.7,1.5) {};
        \node [unmarked={},label=above:{\small4}] (N2) at (2,1.5) {};
        
        \node [unmarked={},label=above:{\small0}] (C0) at (0,3) {};
        \node [unmarked={},label=above:{\small1}] (C1) at (-2,1.5) {};
        \node [unmarked={},label=above:{\small2}] (C2) at (-0.7,1.5) {};
        
    \draw[red](C0) -- (C1);
    \draw[red](C0) -- (C2);
    
    \draw[red](C1) -- (M1);
    \draw[red](C1) -- (M2);
    \draw[red](C1) -- (M3);
    \draw[red](C2) -- (M4);
    \draw[red](C2) -- (M5);
    \draw[red](C2) -- (M6);

    \draw[red](C0) -- (N1);
    \draw[red](C0) -- (N2);
    \draw(N1) -- (M1);
    \draw(N1) -- (M2);
    \draw(N1) -- (M4);
    \draw(N2) -- (M3);
    \draw(N2) -- (M5);
    
    \draw(C2) -- (M1);
    \draw(C2) -- (M3);
    \draw(M1) -- (M2);
    \draw(M3) -- (M4);
    \draw(M5) -- (M6);
    \draw(N1) -- (N2);

    \end{tikzpicture}
\captionsetup{width=0.85\linewidth}
  \caption{A triconed graph. The red edges are the edges of the canonical spanning tree $B_0$.}
\label{fig:exgraph2}  
\end{center}

\vfill\null
\columnbreak

\begin{center}

        \tikzset{
        unmarked/.style args={}{
            draw,fill = black,circle,inner sep = 0pt, minimum size = .30cm,
        }
    }
    
    \begin{tikzpicture}[scale=1]
        \node [unmarked={},label=below:{\small5}] (M1) at (-2.5,0) {};
        \node [unmarked={},label=below:{\small6}] (M2) at (-1.5,0) {};
        \node [unmarked={},label=below:{\small7}] (M3) at (-0.5,0) {};
        \node [unmarked={},label=below:{\small8}] (M4) at (0.5,0) {};
        \node [unmarked={},label=below:{\small9}] (M5) at (1.5,0) {};
        \node [unmarked={},label=below:{\small10}] (M6) at (2.5,0) {};
        \node [unmarked={},label=above:{\small3}] (N1) at (0.7,1.5) {};
        \node [unmarked={},label=above:{\small4}] (N2) at (2,1.5) {};
        \node [unmarked={},label=above:{\small1}] (C1) at (-2,1.5) {};
        \node [unmarked={},label=above:{\small2}] (C2) at (-0.7,1.5) {};

    \draw (N1) -- (M1);
    \draw (N1) -- (M2);
    \draw (N1) -- (M4);
    \draw (N2) -- (M3);
    \draw (N2) -- (M5);
    
    \draw (C2) -- (M1);
    \draw (C2) -- (M3);
    \draw (M1) -- (M2);
    \draw (M3) -- (M4);
    \draw (M5) -- (M6);
    \draw (N1) -- (N2);

    \end{tikzpicture}
    \captionsetup{width=0.85\linewidth}
    
    \caption{The graph $\GG$ obtained from $G$ by removing edges of $B_0$ and vertex $0$.}
    \label{fig:exgraphcone}
    \end{center}

\vfill\null
\columnbreak
\end{multicols}

\end{figure}

\subsection{Trirooted forests}
Given a spanning tree $B$ of a triconed graph $G$, we will show that we can express $B$ as a spanning forest on $\GG$ by taking the edges in $B \cap \GG$ and then adding extra information called \newword{marks} on vertices that encode which edges of $B_0$ were deleted as we go from $G$ to $\GG$. Recall that in $G$, we called the vertices $0,1,2$ to be special and other vertices to be normal. We carry this over to $\GG$ as well, calling the vertices $1,2$ to be special and others to be normal.



For each edge of $B \cap B_0$ excluding $01$ and $02$, the idea is to put a mark on its child. We will elaborate on the details of this process after we define trirooted forests in Definition~\ref{def:trirootedproperties}. The \newword{type} of a mark will be the type of the vertex it is on, as defined in the previous subsection. A \newword{doubly marked vertex} is a vertex of type $2$, that has two marks on it. We think of it as having a usual mark on it of type $2$ and an \newword{extra mark} of type $0$ on the same vertex as well. The purpose of having such object is to encode the information regarding the edge $02$. The type of a tree will be the collection of types of the marks on it. We draw the marks as hollow circles.

For example, take a look at the forest in Figure~\ref{fig:triex3} that has some marks on it. The hollow circles stand for the marks. Vertex $5$ has a mark on it, so the mark here is of type $1$. The vertex $8$ is a doubly marked vertex, and has a mark of type $2$ and a mark of type $0$ on it. The type of the component containing vertex $8$ will be $\{0,1,2\}$.


\begin{definition}
Given a triconed graph $G$, look at the reduced graph $\GG$. If $T$ is a tree on $\GG$ with some marks that satisfy the following properties, we call it a \newword{correctly marked tree} of $\GG$:
\begin{itemize}
    \item there is at least one mark,
    \item special vertices, if contained, are marked,
    \item the types of the marks that appear are all different, 
    \item if there is a doubly marked vertex, there has to be another marked vertex.
\end{itemize}
\end{definition}

Given a tree, let the \newword{degree} of the tree count the number of marks on it (hence if it contains a doubly marked vertex, that contributes two to the degree). Given a component of a forest that has degree $i$, we will simply call it an \newword{i-component} of the forest.


A \newword{blueprint $\P(F)$} of a trirooted forest $F$ will be a graph on vertices $\{0,1,2\}$ with parallel edges allowed coming from the components of $\GG$. For each $2$-component having marks of type $a,b$ we add an edge $ab$ into the blueprint . If there is a $3$-component, we put edges $01$ and $02$ into the blueprint. Beware that unlike $G$, we are allowing parallel edges in the construction of the blueprint.

For example take a look at the forest with marks on it in Figure~\ref{fig:triex3}. Each component is correctly marked. There are no $2$-components and the lone $3$-component has type $\{0,1,2\}$. From this we can construct its blueprint in Figure~\ref{fig:blueprint}: it consists of edges $01$ and $02$.

\begin{definition}
\label{def:trirootedproperties}
Given a triconed graph $G$, a \newword{trirooted forest} (of $\GG$) is a spanning forest $F$ of $\GG$ with some marks on vertices such that the following properties are satisfied:
\begin{itemize}
    \item each component is a correctly marked tree,
    \item its blueprint is acyclic.
\end{itemize}

We will use $\R(\GG)$ to denote the set of trirooted forests of $\GG$.

\end{definition}



Now we will describe how to transform the spanning trees of $G$ to trirooted forests of $\GG$. Suppose $B \in \T(\G)$ is a spanning tree. Remove all edges of $B_0 \cap B$ and when we do so, for all edges being deleted except for $01$ and $02$, mark its child. Mark the vertices $1$ and $2$ regardless of whether $B$ contains $01$ and $02$ or not. If $B$ contained $02$ and there is a path from $2$ to $1$ in $B$ not going through $02$, then pick the first edge along this path that is not in $B_0$, and doubly mark the endpoint of this edge closer to $2$. Notice that this vertex we have chosen, should have been marked already in the first step, so we are essentially adding an extra mark coming from $02$ to that vertex. Every component must have at least one marking since $B$ is connected and spanning. In addition, note that if two marked vertices within the same component have same type, or if two components' image in the blueprint overlap at 2 or more vertices, this implies a cycle in $G$. Hence, this gives us a trirooted forest $F_B$. Let $\phi_1$ be the map that maps $B \in \B(\G)$ to $F_B \in \R(\GG)$.
 
For example, look at the triconed graph in Figure~\ref{fig:triex1}. An example of a spanning tree $B$ of $G$ is given in Figure~\ref{fig:triex2}. For each edge $e$ of $B_0 \cap B$ aside from $02$, we mark its child: from $15$ we mark $5$, from $16$ we mark $6$, from $28$ we mark $8$, and from $2A$ we mark $A$. Now for $02$, look for the path from $2$ to $1$ in $B$. The first vertex along this path going from $2$ to $1$ that is not in $B_0$, turns out to $8$: so we doubly mark it. The blueprint turns out to be Figure~\ref{fig:blueprint} where edges $01$ and $02$ come from the $3$-component.


\begin{figure}[h]
\begin{multicols}{2}
 \begin{center}
 
        \tikzset{
        unmarked/.style args={}{
            draw,fill = black,circle,inner sep = 0pt, minimum size = .30cm,
        }
    }
    
    \begin{tikzpicture}[scale=1]
        \node [unmarked={},label=below:{\small5}] (M1) at (-2.5,0) {};
        \node [unmarked={},label=below:{\small6}] (M2) at (-1.5,0) {};
        \node [unmarked={},label=below:{\small7}] (M3) at (-0.5,0) {};
        \node [unmarked={},label=below:{\small8}] (M4) at (0.5,0) {};
        \node [unmarked={},label=below:{\small9}] (M5) at (1.5,0) {};
        \node [unmarked={},label=below:{\small10}] (M6) at (2.5,0) {};
        \node [unmarked={},label=above:{\small3}] (N1) at (0.5,1.5) {};
        \node [unmarked={},label=above:{\small4}] (N2) at (1.5,1.5) {};
        
        \node [unmarked={},label=above:{\small0}] (C0) at (0,3) {};
        \node [unmarked={},label=above:{\small1}] (C1) at (-1.5,1.5) {};
        \node [unmarked={},label=above:{\small2}] (C2) at (-0.5,1.5) {};
        
    \draw [red](C0) -- (C1);
    \draw [red](C0) -- (C2);
    
    \draw [red](C1) -- (M1);
    \draw [red](C1) -- (M2);
    \draw [red](C1) -- (M3);
    \draw [red](C2) -- (M4);
    \draw [red](C2) -- (M5);
    \draw [red](C2) -- (M6);

    \draw [red](C0) -- (N1);
    \draw [red](C0) -- (N2);
    \draw (N1) -- (M1);
    \draw (N1) -- (M2);
    \draw (N1) -- (M4);
    \draw (N2) -- (M3);
    \draw (N2) -- (M5);
    
    \draw (C2) -- (M1);
    \draw (C2) -- (M3);
    \draw (M1) -- (M2);
    \draw (M3) -- (M4);
    \draw (M5) -- (M6);
    \draw (N1) -- (N2);

    \end{tikzpicture}
       \captionsetup{width=1.0\linewidth}
  \captionof{figure}{An example triconed graph $\G$, with edges of $B_0$ marked in red.}
  \label{fig:triex1}
    
 \end{center}

\vspace{0.5cm}

 \begin{center}
     \tikzset{
        unmarked/.style args={}{
            draw,fill = black,circle,inner sep = 0pt, minimum size = .30cm,
        }
    }

    \tikzset{
        marked/.style args={}{
           draw, circle, fill=none, inner sep=0pt, minimum size=0.30cm, line width=0.04cm
        }
    }
    
    \tikzset{
        marked2/.style args={}{
           draw, double, circle, fill=none, inner sep=0pt, minimum size=0.30cm, line width=0.04cm
        }
    }

    \begin{tikzpicture}[scale=1]
        \node [marked={},label=below:{\small5}] (M1) at (-2.5,0) {};
        \node [marked={},label=below:{\small6}] (M2) at (-1.5,0) {};
        \node [unmarked={},label=below:{\small7}] (M3) at (-0.5,0) {};
        \node [marked2={},label=below:{\small8}] (M4) at (0.5,0) {};
        \node [unmarked={},label=below:{\small9}] (M5) at (1.5,0) {};
        \node [marked={},label=below:{\small10}] (M6) at (2.5,0) {};
        \node [unmarked={},label=above:{\small3}] (N1) at (0.5,1.5) {};
        \node [unmarked={},label=above:{\small4}] (N2) at (1.5,1.5) {};
        
        \node [marked={},label=above:{\small1}] (C1) at (-1.5,1.5) {};
        \node [marked={},label=above:{\small2}] (C2) at (-0.5,1.5) {};
    
    \draw (N1) -- (M2);
    \draw (N1) -- (M4);
    \draw (N2) -- (M3);
    
    \draw (M3) -- (M4);
    \draw (M5) -- (M6);

    \end{tikzpicture}
    \captionsetup{width=1.0\linewidth}
  \captionof{figure}{The result after applying $\phi_1$ to $B$. Marked vertices are denoted by hollow circles, and the doubly marked vertex $8$ is denoted by a double circle.}

  \label{fig:triex3}
 \end{center}

\vfill\null
\columnbreak
 \begin{center}
     \tikzset{
        unmarked/.style args={}{
            draw,fill = black,circle,inner sep = 0pt, minimum size = .30cm,
        }
    }

    \tikzset{
        marked/.style args={}{
           draw, circle, fill=none, inner sep=0pt, minimum size=0.30cm, line width=0.04cm
        }
    }
    
    \tikzset{
        marked2/.style args={}{
           draw, double, circle, fill=none, inner sep=0pt, minimum size=0.30cm, line width=0.04cm
        }
    }

    \begin{tikzpicture}[scale=1]
        \node [unmarked={},label=below:{\small5}] (M1) at (-2.5,0) {};
        \node [unmarked={},label=below:{\small6}] (M2) at (-1.5,0) {};
        \node [unmarked={},label=below:{\small7}] (M3) at (-0.5,0) {};
        \node [unmarked={},label=below:{\small8}] (M4) at (0.5,0) {};
        \node [unmarked={},label=below:{\small9}] (M5) at (1.5,0) {};
        \node [unmarked={},label=below:{\small10}] (M6) at (2.5,0) {};
        \node [unmarked={},label=above:{\small3}] (N1) at (0.5,1.5) {};
        \node [unmarked={},label=above:{\small4}] (N2) at (1.5,1.5) {};
        
        \node [unmarked={},label=above:{\small0}] (C0) at (0,3) {};
        \node [unmarked={},label=above:{\small1}] (C1) at (-1.5,1.5) {};
        \node [unmarked={},label=above:{\small2}] (C2) at (-0.5,1.5) {};
    
    \draw [red](C1) -- (M1);
    \draw [red](C1) -- (M2);
    \draw [red](C2) -- (M4);
    \draw [red](C2) -- (M6);

    \draw [red](C0) -- (C2);
    \draw (N1) -- (M2);
    \draw (N1) -- (M4);
    \draw (N2) -- (M3);
    
    \draw (M3) -- (M4);
    \draw (M5) -- (M6);

    \end{tikzpicture}
    \captionsetup{width=1.0\linewidth}
  \captionof{figure}{An example spanning tree $B$ of $\G$, with its cone edges in red.}

  \label{fig:triex2}
 \end{center}

 \begin{center}
     \tikzset{
        unmarked/.style args={}{
            draw,fill = black,circle,inner sep = 0pt, minimum size = .30cm,
        }
    }

    \begin{tikzpicture}[scale=1]
        
        \node [unmarked={},label=above:{\small0}] (C0) at (0,3) {};
        \node [unmarked={},label=above:{\small1}] (C1) at (-1.5,1.5) {};
        \node [unmarked={},label=above:{\small2}] (C2) at (-0.5,1.5) {};

    \draw [black](C0) -- (C1);
    \draw [black](C0) -- (C2);

    \end{tikzpicture}
    \captionsetup{width=1.0\linewidth}
  \captionof{figure}{The blueprint of $B$. $B$ contains a 3-component, so the edges 01 and 02 are present.}

  \label{fig:blueprint}
 \end{center}

\end{multicols}

\end{figure}

As we go from $B$ to $\phi_1(B)$, there is a bijection between the edges of $B_0 \cap B$ excluding $01$ and $02$, and the marked vertices in $\phi_1(B)$. For edges $01$ and $02$, their existence is encoded in the blueprint as we can see in the following lemmas:

\begin{lemma}
\label{lem:01blueprint}
Let $B$ be a spanning tree of a triconed graph $G$. We have $01$ as an edge inside $B$ if and only if $0$ and $1$ are not connected in the blueprint. 
\end{lemma}
\begin{proof}
If $01$ is in $B$, then there cannot exist a component in $\phi_1(B)$ that has type set $\supseteq \{0,1\}$, since otherwise we have a cycle in $B$. Hence $01$ cannot exist inside the blueprint. If we have a component of type $\{1,2\}$ and a component of type $\{0,2\}$ in $\phi_1(B)$, we again get a cycle inside $B$ using the edge $01$. Therefore $0$ and $1$ will not be connected in the blueprint.

If we don't have $01$ inside $B$, there still is a path from $0$ to $1$ inside $B$ since it is a spanning tree. If there exists a component of type $\supseteq \{0,1\}$ in $\phi_1(B)$ then we have $01$ inside the blueprint. If there is no such component, we must have some component that connects $0$ to $2$ and some component that connects $1$ to $2$ in $\phi_1(B)$, which means having $02$ and $12$ inside the blueprint.
\end{proof}

\begin{lemma}
\label{lem:02blueprint}
Let $B$ be a spanning tree of a triconed graph $G$ that doesn't contain $01$ but does contain $02$. There is a path within $B$ from $1$ to $2$ that doesn't go through $02$ if and only if $1$ and $2$ are connected in the blueprint.
\end{lemma}
\begin{proof}
First notice that since $02$ is an edge of $B$, we cannot have $0$ and $2$ connected inside the blueprint. Hence existence of a path from $1$ to $2$ that doesn't go through $02$ in $B$ is equivalent to the existence of a component of type $\{1,2\}$ in $\phi_1(B)$.
\end{proof}

For example take a look at the blueprint in Figure~\ref{fig:blueprint}. From Lemma~\ref{lem:01blueprint} tells us that $01$ is not present in the corresponding spanning tree. From Lemma~\ref{lem:02blueprint} tells us that if the corresponding spanning tree contains $02$, then there is a path from $1$ to $2$ that doesn't go through $02$.  We can verify that these statements are true in Figure~\ref{fig:triex2}.

\begin{proposition}
\label{trirooted}
The map $\phi_1: \T(\G) \rightarrow \R(\GG)$ described above is a bijection. 
\end{proposition}

\begin{proof}
All we need to do is describe how to recover the spanning tree from a trirooted forest. Given a trirooted forest $F$, take all normal marked vertices and add their cone edges accordingly. If there is a doubly-marked vertex, add in the edge $02$. Using Lemma~\ref{lem:01blueprint}, if $0$ and $1$ were not connected in the blueprint add in $01$ to our construction. If we still do not have a tree, then add in $02$.

\end{proof}

For example, using the algorithm described in the proof above we can go from Figure~\ref{fig:triex3} back to Figure~\ref{fig:triex2}. For the first step, we take all normal marks $5,6,8,A$ and add their cone edges $15,16,28,2A$. Now since there is a doubly-marked vertex, we add in the edge $02$. Since the blueprint was spanning, we do not need to add the edge $01$. At this point we are done reconstructing the tree of Figure~\ref{fig:triex2}.


Now that we have established that $\phi_1$ is a bijection between spanning trees of $G$ and trirooted forests of $\GG$, we will show that it carries over the information regarding passivity nicely.

\begin{lemma}
\label{lem:edgepassive}
Let $B$ be a spanning tree of a triconed graph $G$. The cone edge of a marked vertex in a $1$-component where the vertex is also the smallest vertex of the component is active. The edge $01$ is active as well. The edge $02$ is active when there is no doubly marked vertex. All other edges are passive.
\end{lemma}

\begin{proof}
To start, notice that all edges of $B \setminus B_0$ are passive thanks to the matroid exchange property: we can exchange any edge in $B \setminus B_0$ with some (smaller) edge of $B_0$.

The edge $01$ can never be passive since it is the smallest edge. The edge $02$ is passive if and only if there is a doubly marked vertex: if there is a doubly marked vertex, we can replace $02$ with $01$.


Now we will describe exactly what edges of $B \cap B_0 \setminus \{01,02\}$ are passive. Pick an edge $xv$ where $v$ is a marked vertex. If $v$ is in a $3$-component it is fairly simple: if $v$ is of type $2$, we can replace the edge $xv$ with $02$, otherwise we can replace it with $01$. If $v$ is in a $2$-component of type $\{0,1\}$ we can replace the edge $xv$ with $01$, if of type $\{0,2\}$ we replace with $02$. In the case $v$ is in a $2$-component of type $\{1,2\}$ we need to look at the blueprint. We have an edge $12$ in the blueprint coming from this component. If we have the edge $02$ in the blueprint, we replace the edge $xv$ in question with $01$ and if not, we replace it with $02$.

Now the only remaining case is when $v$ is in a $1$-component. If $v$ is the smallest vertex of the component, the edge is active. If not, we can replace the edge $xv$ with the cone edge of the smallest vertex of the component. 
\end{proof}

For example take a look at the tree $B$ in Figure~\ref{fig:exampletree1}. This is a spanning tree of the triconed graph $G$ from Figure~\ref{fig:triex1}. Edges in $B \setminus B_0$, the edges $36,38,47,78,9A$, are passive. Among the edges of $B \cap B_0$, the edges $16,02,28,2A$ are passive and the edge $15$ is active. The edge $15$ is active since its corresponding mark at vertex $5$, is on the smallest vertex of that component which only contains $5$. The edges $16,28,2A$ are passive since their corresponding marks on vertices $6,8,A$ are in components that have at least two marks. Finally the edge $02$ is passive since there is a doubly marked vertex.

Given a mark on a vertex $v$, recall that it corresponds to a unique edge of $B_0$: namely the edge that has the vertex $v$ as its child (the exception is the extra mark on the doubly marked vertex, which corresponds to the edge $02$). We say that a mark is a \newword{passive mark} if the associated edge is passive, and an \newword{active mark} if the associated edge is active. From the above lemma we get the following corollary:

\begin{corollary}
\label{cor:maxpassivespanning}
Let $B$ be a spanning tree of a triconed graph $G$ and $\phi_1(B)$ be its corresponding trirooted forest in $\GG$. All edges of $B$ are passive if and only if $\P(\phi_1(B))$ is spanning and all normal marks of $\phi_1(B)$ are passive.
\end{corollary}
\begin{proof}
Aside from $01$ and $02$, any edge of $B_0$ is passive if and only if the corresponding mark is a passive mark. The edge $01$, if it exists in a tree is always active. The existence of $01$ is equivalent to $0$ and $1$ not being connected in the blueprint from Lemma~\ref{lem:01blueprint}. The edge $02$ is passive if and only if $01$ is missing and there is a doubly marked vertex, which results in a spanning blueprint: if the blueprint is not spanning, then that means that an active $01$ or $02$ is present. Conversely, if the blueprint is spanning, then neither an active $01$ nor an active $02$ are present.

\end{proof}

\begin{figure}[h]
\begin{multicols}{2}

 \begin{center}
     \tikzset{
        unmarked/.style args={}{
            draw,fill = black,circle,inner sep = 0pt, minimum size = .30cm,
        }
    }

    \tikzset{
        marked/.style args={}{
           draw, circle, fill=none, inner sep=0pt, minimum size=0.30cm, line width=0.04cm
        }
    }
    
    \tikzset{
        marked2/.style args={}{
           draw, double, circle, fill=none, inner sep=0pt, minimum size=0.30cm, line width=0.04cm
        }
    }

    \begin{tikzpicture}[scale=1]
        
        \node [unmarked={},label=below:{\small5}] (M1) at (-2.5,0) {};
        
        \node [unmarked={}, label=below:{\small6}] (M2) at (-1.5,0) {};
        \node [unmarked={}, label=below:{\small7}] (M3) at (-0.5,0) {};
        \node [unmarked={},   label=below:{\small8}] (M4) at (0.5,0) {};
        \node [unmarked={}, label=below:{\small9}] (M5) at (1.5,0) {};
        \node [unmarked={},   label=below:{\small10}] (M6) at (2.5,0) {};
        \node [unmarked={}, label=above:{\small3}] (N1) at (0.5,1.5) {};
        \node [unmarked={},   label=above:{\small4}] (N2) at (1.5,1.5) {};
        
        \node [unmarked={},label=above:{\small0}] (C0) at (0,3) {};
        \node [unmarked={},label=above:{\small1}] (C1) at (-1.5,1.5) {};
        \node [unmarked={},label=above:{\small2}] (C2) at (-0.5,1.5) {};
    
    \draw (C1) -- (M1);
    \draw [red](C1) -- (M2);
    \draw [red](C2) -- (M4);
    \draw [red](C2) -- (M6);

    \draw [red](C0) -- (C2);
    \draw [red](N1) -- (M2);
    \draw [red](N1) -- (M4);
    \draw [red](N2) -- (M3);
    
    \draw [red](M3) -- (M4);
    \draw [red](M5) -- (M6);

    \end{tikzpicture}
    \captionsetup{width=1.0\linewidth}
  \captionof{figure}{An example spanning tree $B$ of $\G$, now with its passive edges in red. There are $10$ passive edges, $4$ of which are in $B\cap B_0$.}

  \label{fig:exampletree1}
 \end{center}
 
\vfill\null
\columnbreak

 \begin{center}
     \tikzset{
        unmarked/.style args={}{
            draw,fill = black,circle,inner sep = 0pt, minimum size = .30cm,
        }
    }

    \tikzset{
        marked/.style args={}{
           draw, circle, fill=none, inner sep=0pt, minimum size=0.30cm, line width=0.04cm
        }
    }
    
    \tikzset{
        marked2/.style args={}{
           draw, double, circle, fill=none, inner sep=0pt, minimum size=0.30cm, line width=0.04cm
        }
    }

    \begin{tikzpicture}[scale=1]
        \node [marked={},label=below:{\small5}] (M1) at (-2.5,0) {};
        \node [marked={},label=below:{\small6}] (M2) at (-1.5,0) {};
        \node [unmarked={},label=below:{\small7}] (M3) at (-0.5,0) {};
        \node [marked2={},label=below:{\small8}] (M4) at (0.5,0) {};
        \node [unmarked={},label=below:{\small9}] (M5) at (1.5,0) {};
        \node [marked={},label=below:{\small10}] (M6) at (2.5,0) {};
        \node [unmarked={},label=above:{\small3}] (N1) at (0.5,1.5) {};
        \node [unmarked={},label=above:{\small4}] (N2) at (1.5,1.5) {};
        
        \node [marked={},label=above:{\small1}] (C1) at (-1.5,1.5) {};
        \node [marked={},label=above:{\small2}] (C2) at (-0.5,1.5) {};
    
    \draw (N1) -- (M2);
    \draw (N1) -- (M4);
    \draw (N2) -- (M3);
    
    \draw (M3) -- (M4);
    \draw (M5) -- (M6);

    \end{tikzpicture}
    \captionsetup{width=1.0\linewidth}
  \captionof{figure}{The result after applying $\phi_1$ to $B$. There are $7$ markings in total, and $1$ marked vertex in a $1$-component which is the smallest vertex of the component.} 

  \label{fig:exampledoubletree1}
 \end{center}

\end{multicols}

\end{figure}


\begin{figure}[h]
\begin{multicols}{3}
 \begin{center}
 
        \tikzset{
        unmarked/.style args={}{
            draw,fill = black,circle,inner sep = 0pt, minimum size = .30cm,
        }
    }

    \tikzset{
        marked/.style args={}{
           draw, circle, fill=none, inner sep=0pt, minimum size=0.30cm, line width=0.04cm
        }
    }
    
    \tikzset{
        marked2/.style args={}{
           draw, double, circle, fill=none, inner sep=0pt, minimum size=0.30cm, line width=0.04cm
        }
    }

    \begin{tikzpicture}[scale=1]
        
        \node [marked={},label=below:{\small5}] (M1) at (-2.5,0) {};
        
        \node [marked={}, label=below:{\small6}] (M2) at (-1.5,0) {};
        \node [unmarked={}, label=below:{\small7}] (M3) at (-0.5,0) {};
        \node [marked={},   label=below:{\small8}] (M4) at (0.5,0) {};
        \node [unmarked={}, label=below:{\small9}] (M5) at (1.5,0) {};
        \node [marked={},   label=below:{\small10}] (M6) at (2.5,0) {};
        \node [unmarked={}, label=above:{\small3}] (N1) at (0.5,1.5) {};
        \node [marked={},   label=above:{\small4}] (N2) at (1.5,1.5) {};
        
        \node [unmarked={},label=above:{\small0}] (C0) at (0,3) {};
        \node [marked={},label=above:{\small1}] (C1) at (-1.5,1.5) {};
        \node [marked={},label=above:{\small2}] (C2) at (-0.5,1.5) {};
    
    \draw [red](C1) -- (M1);
    \draw [red](C1) -- (M2);
    \draw [red](C2) -- (M4);
    \draw [red](C2) -- (M6);

    \draw [red](C0) -- (N2);
    \draw (N1) -- (M2);
    \draw (N1) -- (M4);
    \draw (N2) -- (M3);
    
    \draw (M3) -- (M4);
    \draw (M5) -- (M6);

    \end{tikzpicture}
    \captionsetup{width=1.0\linewidth}
  \captionof{figure}{
     Another spanning tree $B$ of $\G$, with its cone edges in red and its marked vertices denoted by hollow circles.}
  \label{fig:spantree}
 \end{center}

\vfill\null
\columnbreak
 \begin{center}
 
        \tikzset{
        unmarked/.style args={}{
            draw,fill = black,circle,inner sep = 0pt, minimum size = .30cm,
        }
    }

    \tikzset{
        marked/.style args={}{
           draw, circle, fill=none, inner sep=0pt, minimum size=0.30cm, line width=0.04cm
        }
    }
    
    \tikzset{
        marked2/.style args={}{
           draw, double, circle, fill=none, inner sep=0pt, minimum size=0.30cm, line width=0.04cm
        }
    }

    \begin{tikzpicture}[scale=1]
        
        \node [marked={},label=below:{\small5}] (M1) at (-2.5,0) {};
        
        \node [marked={}, label=below:{\small6}] (M2) at (-1.5,0) {};
        \node [unmarked={}, label=below:{\small7}] (M3) at (-0.5,0) {};
        \node [marked={},   label=below:{\small8}] (M4) at (0.5,0) {};
        \node [unmarked={}, label=below:{\small9}] (M5) at (1.5,0) {};
        \node [marked={},   label=below:{\small10}] (M6) at (2.5,0) {};
        \node [unmarked={}, label=above:{\small3}] (N1) at (0.5,1.5) {};
        \node [marked={},   label=above:{\small4}] (N2) at (1.5,1.5) {};
        
        \node [marked={},label=above:{\small1}] (C1) at (-1.5,1.5) {};
        \node [marked={},label=above:{\small2}] (C2) at (-0.5,1.5) {};
    
    \draw (N1) -- (M2);
    \draw (N1) -- (M4);
    \draw (N2) -- (M3);
    
    \draw (M3) -- (M4);
    \draw (M5) -- (M6);

    \end{tikzpicture}
       \captionsetup{width=1.0\linewidth}
  \captionof{figure}{The associated trirooted forest.}
  \label{fig:trirooted}
 \end{center}
 
 \vfill\null
\columnbreak
 \begin{center}
 
        \tikzset{
        unmarked/.style args={}{
            draw,fill = black,circle,inner sep = 0pt, minimum size = .30cm,
        }
    }

    \tikzset{
        single/.style args={}{
            draw,line width=0.5mm,black
        }
    }
    \tikzset{
        double/.style args={}{
            preaction={draw,line width=1.25mm,black},draw,line width=.25mm,white
        }
    }
    \tikzset{
        triple/.style args={}{
            preaction={preaction={draw,line width=2mm,black},draw,line width=1mm,white},draw,line 
            width=.5mm,black
        }
    }  
    \tikzset{
        quadruple/.style args={}{
            preaction={preaction={preaction={draw,line width=2.75mm,black},draw,line width=1.75mm,white},draw,line width=1.25mm,black},draw,line width=.25mm,white
        }
    }

    \begin{tikzpicture}[scale=1]
        
        \node [unmarked={},label=below:{\small5}] (M1) at (-2.5,0) {};
        
        \node [unmarked={}, label=below:{\small6}] (M2) at (-1.5,0) {};
        \node [unmarked={}, label=below:{\small7}] (M3) at (-0.5,0) {};
        \node [unmarked={}, label=below:{\small8}] (M4) at (0.5,0) {};
        \node [unmarked={}, label=below:{\small9}] (M5) at (1.5,0) {};
        \node [unmarked={}, label=below:{\small10}] (M6) at (2.5,0) {};
        \node [unmarked={}, label=above:{\small3}] (N1) at (0.5,1.5) {};
        \node [unmarked={}, label=above:{\small4}] (N2) at (1.5,1.5) {};
        
        \node [unmarked={},label=above:{\small1}] (C1) at (-1.5,1.5) {};
        \node [unmarked={},label=above:{\small2}] (C2) at (-0.5,1.5) {};
    
    \draw [double={}](N1) -- (M2);
    \draw [single={}](N1) -- (M4);
    \draw [triple={}](N2) -- (M3);
    
    \draw [single={}](M3) -- (M4);
    \draw [double={}](M5) -- (M6);

    \end{tikzpicture}
      \captionsetup{width=1.0\linewidth}
  \captionof{figure}{The associated 3-weighted forest.}
    \label{fig:3weighted}
  \end{center}

\end{multicols}

\end{figure}


\subsection{$3$-weighted-forests} 
Recall that our goal is to convert a spanning tree to a monomial while carrying over the information of number of passive edges into the degree. Hence we want to convert the information of marks (which were placed on vertices) onto weights on the edges of the graph, which is what we will do in this subsection.

We will use the same notation as \cite{Kook2007} and \cite{biconed}: let $\G$ be a triconed graph and suppose $F$ is any set of edges of $\GG$. Let $\omega_F:F \rightarrow {\mathbb N}_{\geq 1}$ denote a positive integer weighting on these set of edges. We will use $\omega$ if the context is clear and refer to the pair $(F,\omega)$ as a \newword{weighted collection of edges}. We say $e \in F$ has \newword{weight n} if $\omega_F(e) = n$. 

\begin{remark}
Notice that a weighted collection of edges naturally gives us a monomial. We can think of the edges of $\GG$ as variables and the weights on the edges as the powers of the variables. For example the monomial associated to the weighted forest in Figure~\ref{fig:3weighted} would be $x_{36}^2x_{38}x_{78}x_{47}^3x_{9A}.$ Hence we will go between a weighted forest and its corresponding monomial freely: for example, whenever we say that a weighted forest of $\GG$ divides another weighted forest of $\GG$, we mean division between the corresponding monomials.
\end{remark}

Given a weighted collection of edges $(F, \omega_F)$, let us define a \newword{weighted component} (or just \newword{component} if it is clear from the context) as a connected component $S \subset F$, along with its weighting $\omega_S$ which is simply equal to $\omega_F$ restricted to $S$. Then, we calculate the \newword{excess weight} of a weighted component $S$ as
\[\delta_{1 \in S} + \delta_{2 \in S} + \sum_{e \in S} (\omega(e) - 1),\] 
where $\delta_{a \in S}$ equals $1$ if $a \in S$ and equals $0$ otherwise. We call a component $S$ to be \newword{$k$-xweighted} if its excess weight is $k$ (note that the excess weight of a component is always nonnegative).

For example, let us look at the 3-weighted forest in figure \ref{fig:3weighted}. The excess weight of the component containing vertices $3,4,6,7,$ and $8$ is $3$, and hence is $2$-xweighted. The excess weight of the component containing the vertex $1$, and the component containing the vertex $2$, is $1$. All other components have excess weight $0$.

\begin{remark}
Excess weight will turn out to be same as the degree of a trirooted component. Weight will turn out to be same as the degree of the monomial coming from that component.
\end{remark}

The main idea is to convert a component of a trirooted forest to a component of a $3$-weighted forest, while keeping the edge set. Hence given a component, the information stored on the marks of the vertices of that component will be stored on the weights of the edges of that component. We start by defining the analogue of the correctly marked trees:

\begin{definition}
Let $(T,w)$ be a weighted tree of $\GG$. We say that it is a \newword{correctly weighted tree} if it satisfies the following properties:
\begin{itemize}
    \item [$(C1)$] Has excess weight $\leq 3$.
    \item [$(C2)$] Let us define a \newword{weighted object} as either an edge of $2$ or more weight, or the vertices $1$ or $2$. Then, if a component has at least $2$ excess weight, any pair of weighted objects has to be \newword{compatible}: the (unique) shortest path containing both of the weighted objects must have endpoints of different type. If an edge has weight at least $3$, we consider the edge to contain distinct weighted objects based on the same edge.
    \item[$(C3)$] The maximum weight of an edge is the sum of the heights of its endpoints. If we have an edge that has its weight equal to that maximal value, we call it a \newword{maximal edge}.
    \item[$(C4)$] If there is a maximal edge, all other edges are at most $2$-weighted.    
\end{itemize}
\end{definition}

We start out by describing a map $\phi_2$ that sends a correctly marked tree of $\GG$ to a correctly weighted tree of $\GG$. By default we set the weight of all edges as $1$. Now each mark will potentially \newword{empower} an associated edge: meaning that it will contribute an extra weight to that edge. Moreover we want each passive mark to empower an edge, and active marks do not. So the number of passive marks we started out with will equal the total number of empowerments done. 

First in the case we are looking at a degree $1$-component in a trirooted forest, let $v$ be the marked vertex and $w$ be the smallest vertex of the component. If $v=w$, that marked vertex is not passive, hence this mark will not empower any edge and we keep all the weights of the component to be $1$. Otherwise, we look at the unique path within the component between between $v$ and $w$, take the edge closest to $v$ within that path and empower that edge with the mark at $v$, setting its weight to $2$.

If we have a degree $2$-component in a trirooted forest, let $v$ and $w$ be the marked vertices. Since doubly marked vertices can only occur in degree $3$-components, we may safely assume that $v \not = w$. Similar to above, we look at the unique path between $v$ and $w$ within the component and associate the mark of $v$ to the edge of the path closest to $v$. Similarly we associate the mark of $w$ to the edge of that path closest to $w$. Now if $v$ is normal then we empower the associated edge, then do the same for $w$. Notice that regardless of whether $v$ is special or not, the path between the weighted objects that forms from this process stays the same as the path between the marked vertices.

The remaining case is when we are looking at a $3$-component in a trirooted forest. We first go over the case when we do not have a doubly marked vertex: let $v,w,z$ be the marked vertices that are all different. If the unique minimal tree that contains $v,w,z$ in the component has all $v,w,z$ as its leaf, we do the same procedure as above: for $v$ (respectively for $w$ and then $z$) associate the edge within the tree adjacent to $v$ (respectively for $w,z$), then empower the edge if $v$ is normal (respectively for $w,z$). If not, then change the labels so that $z$ lies among the unique path $P$ between $v$ and $w$ within the component, and that $\type{v} < \type{w}$ (the case when we have $z=w$, that is when we have a doubly marked vertex at $w$ also falls into this case). From now we will describe this special case as the mark at $z$ being \newword{sandwiched} by the marks at $v$ and $w$. We pick the associated edges of $v$ and $w$ as before: pick the edge adjacent to the vertex within the path. For $z$, within the path from $z$ to $v$ starting from $z$, look for the first edge $e$ such that its endpoint closer to $v$ is compatible with $w$ and its endpoint closer to $w$ is compatible with $v$: whenever there is such edge we will say that $e$ is compatible with $P$. We explain why we can always find such $e$ in the lemma below:

\begin{lemma}
\label{lem:eexist}
In the process above, we can always find an edge compatible with the path $P$ as we move from $z$ to $v$. Moreover, we can find such edge before we reach a vertex that has same type as $z$.
\end{lemma}
\begin{proof}
Assume for sake of contradiction, we cannot find such edge. Starting from the edge that contains $z$, since $\type{z}\neq\type{v}$ it must be incompatible with $w$, so its other endpoint has to be same type as $w$. Now whenever we have an edge where the endpoint closer to $z$ has type $w$, its other endpoint has to be type $w$ again. This leads to all vertices from on the way from $z$ to $v$, except for $z$, have type same as $w$, which contradicts the fact that $\type{v},\type{z} \not = \type{w}$. 
\end{proof}

Hence the above procedure is well defined.

\begin{remark}
\label{rem:markempower}
From our construction, except when we are looking at at the case when we have a sandwiched mark, a mark always empowers the edge adjacent to it. Also each empowerment comes from a passive mark: the excess weight of a component equals the number of passive marks it had before.
\end{remark}

\begin{lemma}
Let $T$ be a correctly marked tree of $\GG$. Then its image $\phi_2(T)$ is a correctly weighted forest of $\GG$. 
\end{lemma}

\begin{proof}
Since each marked vertex either provides $0$ or $1$ excess weight, we get $(C1)$ from the fact that each component contains at most $3$ marks. The condition $(C2)$ follows from the construction: we picked a path or a tree and made sure the endpoints were compatible. To see that condition $(C3)$ and $(C4)$ is satisfied we do a case-by-case analysis. 

First start with the case when we have an edge $e$ where both endpoints are of height $1$. Assume for the sake of contradiction that $w(e) \geq 3$. From Remark~\ref{rem:markempower}, at least one endpoint is a passive mark of type $0$ and call it $v$. Since both endpoints cannot be passive (otherwise they are both of type $0$), the remaining empowerment has to come from a sandwiched mark. Let $z$ be the sandwiched mark that empowers $e$ and $w$ be the remaining mark not of type $0$. Since $e$ is the first edge within the path from $z$ to $v$ that is compatible with the path $vw$, all vertices from $z$ to $v$ except $z$ and $v$ should have same type as $w$. But this means the endpoint of $e$ that is not $v$ has to be a special vertex and $w$ is that vertex. Since $z$ has to be sandwiched between $v$ and $w$, we get that $z=w$ which contradicts that we cannot have a doubly marked vertex without the extra mark being type $0$. Hence $(C3)$ holds in this case. Now when $e$ is a maximal edge, in order for us to have another edge $f$ that is $3$-weighted in the component, we need to have a sandwiched mark: there are three marks $v,z,w$ on a path with $\type{v} < \type{w}$ where the mark on $w$ empowers $e$ and marks at $v,z$ empowers $f$. Since both endpoints of $e$ are of height $1$, we have $w$ to be of type $0$, which contradicts $\type{v} < \type{w}$.



In the case when $e$ has both endpoints of height $2$, its weight being bounded by $4$ is basically $(C1)$ which we have shown. Condition $(C4)$ holds automatically as well from $(C1)$. So the remaining case to consider is when one endpoint of $e$ is of height $1$ and the other being height $2$. We name them $a$ and $b$ respectively. Assume for the sake of contradiction that $w(e)=4$. For this to happen, the component containing $e$ must have $3$ passive marks, meaning no vertex $1$ or $2$ present. Moreover, since $a$ is of type $0$, there is no doubly marked vertex. But then we cannot have three passive marks all empowering $e$, which leads to a contradiction. This finishes the proof of $(C3)$ in this case. Condition $(C4)$ in this case holds automatically from $(C1)$. 
\end{proof}

For example take a look at the trirooted forest in Figure~\ref{fig:trirooted}. First start with the trivial component which only contains the vertex $5$. Its image under $\phi_2$ is again going to be a single vertex $5$. Next let us look at the tree consisting of vertices $9$ and $A$. This has degree $1$ and the marked at $A$ is not on the smallest vertex of the component. This is a passive mark and will empower the edge $9A$. 
Hence the image of this tree under $\phi_2$ will be the edge $9A$ having weight $2$. Finally look at the component consisting of vertices $3,4,6,7,8$. This is a degree $3$ component, with the mark at $8$ being on the path between mark at $6$ and $4$. So the mark at $6$ empowers $36$, and the mark at $8$ empowers $47$ (the edges lying at the end of the path between $6$ and $8$). To see which edge the mark at $8$ empowers, we look at the segment from $8$ to $4$ (here $4$ is chosen instead of $6$ since $4<6$), then look for the first edge such that the endpoints are compatible with marks at $4,6$ respectively. Edge $78$ does not qualify since $\type{7} = \type{6} = 1$. Edge $47$ qualifies since $\type{4} \not = \type{6}$ and $\type{8} \not = \type{4}$. Hence its image under $\phi_2$ will be the same underlying tree with edge $36$ empowered once and $47$ empowered twice. 

\begin{proposition}
The map $\phi_2$ is a bijection between correctly marked trees of $\GG$ and correctly weighted trees of $\GG$.
\end{proposition}
\begin{proof}

Here we describe that inverse map to $\phi_2$: a map that sends a correctly weighted tree to a correctly marked tree. We do a case-by-case analysis based on the number of weighted objects inside the component which we denote as $c$. In the case $c=0$, simply set the smallest vertex in the component as the sole marked vertex. In the case $c=1$, if the weighted object is a special vertex, mark that vertex. If the weighted object is an empowerment on an edge $e$, look at the shortest unique path between the smallest vertex $v$ of the component and $e$, then put a mark on the endpoint $\not = v$. When $c=2$, take the shortest path containing the objects, put marks on the endpoints. Finally in the case $c=3$, take the smallest tree containing these objects. If the empowered objects each contain different leaves of this tree (so there are three leaves), simply mark those leaves. If not then this is a path between vertices $a < b$ with an empowered edge in the middle. From this edge in the middle, move towards $b$ and look for the first vertex that has type different from $a$ and $b$, then mark that vertex and mark the vertices $a$ and $b$ as well.
\end{proof}

Now given a trirooted forest, we define $\phi_2$ on it by applying $\phi_2$ on each of its component (a tree) and taking the union of the resulting weighted trees to get a weighted forest. Given a weighted forest $(F,w)$ where each component is a correctly weighted tree, we define its blueprint $\P(F,w)$ as the blueprint of its preimage under $\phi_2$.


\begin{definition}
\label{def:3erf}
For a triconed graph $\G$, a \newword{3-weighted forest} is a weighted set of edges $(F, \omega)$ with $F \subset E(\GG)$ satisfying:
\begin{itemize}
    \item $F$ induces a forest in $\GG$. Each component is a correctly weighted tree.
    \item The blueprint $\P(F,w)$ is acyclic.
\end{itemize}
\end{definition}

For example take a look at the weighted forest $(F,w)$ in Figure~\ref{fig:3weighted}. We have two nontrivial components, each of them being a correctly weighted tree. To find the blueprint, we look at the preimage of each component under the map $\phi_2$, to get the two nontrivial components in Figure~\ref{fig:trirooted}. Hence the blueprint turns out to be edges $01$ and $02$. It being acyclic implies that what we have in Figure~\ref{fig:3weighted} is indeed a $3$-weighted forest.


\section{Pure multicomplex}

In the previous section we constructed a map $\phi_1$ that sends a spanning tree of a triconed graph $G$ into a trirooted forest, and a map $\phi_2$ that sends a trirooted forest into a $3$-weighted forest which can be identified with a monomial supported on a forest. In this section we study the properties of the image set of $\phi := \phi_2 \circ \phi_1$ and show that we get a pure order ideal.

\begin{proposition}
The image of $\phi$ is a multicomplex.

\end{proposition}
\begin{proof}
Let $T$ be any trirooted forest and $\phi_2(T)$ be the resulting $3$-weighted forest. We want to show that for any edge $e$ that appears in $\phi_2(T)$ (hence in $T$), that $\phi_2(T)/e$ is again a $3$-weighted forest. Let $C$ stand for the component containing $e$. We first need to show that $\phi_2(C)/e$ is either a correctly weighted tree (when $w(e) \geq 2$) or becomes two correctly weighted trees (when $w_e = 1$). The conditions $(C1),(C3),(C4)$ again hold pretty much from definition, so we only have to check for $(C2)$. When $w_e$ was at least $2$, the compatibility again follows. When $w_e$ was $1$, so we are deleting the edge $e$, the set of weighted objects gets partitioned into two disjoint components, so compatibility still holds.

Now we have to check that the blueprint is still acyclic. Unless in the case we have a sandwiched mark, the resulting blueprint is a subgraph of the original blueprint of the component. In the case there is a sandwiched mark, which only appears in a degree $3$ component, recall that there are no other edges in the blueprint of $T$. The new blueprint will be a single edge, hence is automatically acyclic.
\end{proof}

For example take a look at the $3$-weighted forest from Figure~\ref{fig:3weighted}. Thinking of this as a monomial, we want to show that this monomial divided by $x_{36}$ is again a $3$-weighted forest: it is drawn in Figure~\ref{fig:3weighteddivided}. First notice that as we remove the empowerment of $36$ from the component, what we get as a result is still a correctly weighted tree. The preimage of this new tree is in the component consisting of vertices $3,4,6,7,8$ in Figure~\ref{fig:trirooteddivided}. The blueprint of this new component is just the edge $02$, so the blueprint of Figure~\ref{fig:3weighteddivided} is acyclic. Hence the weighted forest in Figure~\ref{fig:3weighteddivided} is indeed a $3$-weighted forest. The corresponding spanning tree of $G$ is drawn in Figure~\ref{fig:spantreedivided}.

\begin{figure}[h]
\begin{multicols}{3}
 \begin{center}
 
        \tikzset{
        unmarked/.style args={}{
            draw,fill = black,circle,inner sep = 0pt, minimum size = .30cm,
        }
    }

    \tikzset{
        single/.style args={}{
            draw,line width=0.5mm,black
        }
    }
    \tikzset{
        double/.style args={}{
            preaction={draw,line width=1.25mm,black},draw,line width=.25mm,white
        }
    }
    \tikzset{
        triple/.style args={}{
            preaction={preaction={draw,line width=2mm,black},draw,line width=1mm,white},draw,line 
            width=.5mm,black
        }
    }  
    \tikzset{
        quadruple/.style args={}{
            preaction={preaction={preaction={draw,line width=2.75mm,black},draw,line width=1.75mm,white},draw,line width=1.25mm,black},draw,line width=.25mm,white
        }
    }

    \begin{tikzpicture}[scale=1]
        
        \node [unmarked={},label=below:{\small5}] (M1) at (-2.5,0) {};
        
        \node [unmarked={}, label=below:{\small6}] (M2) at (-1.5,0) {};
        \node [unmarked={}, label=below:{\small7}] (M3) at (-0.5,0) {};
        \node [unmarked={}, label=below:{\small8}] (M4) at (0.5,0) {};
        \node [unmarked={}, label=below:{\small9}] (M5) at (1.5,0) {};
        \node [unmarked={}, label=below:{\small10}] (M6) at (2.5,0) {};
        \node [unmarked={}, label=above:{\small3}] (N1) at (0.5,1.5) {};
        \node [unmarked={}, label=above:{\small4}] (N2) at (1.5,1.5) {};
        
        \node [unmarked={},label=above:{\small1}] (C1) at (-1.5,1.5) {};
        \node [unmarked={},label=above:{\small2}] (C2) at (-0.5,1.5) {};
    
    \draw [single={}](N1) -- (M2);
    \draw [single={}](N1) -- (M4);
    \draw [triple={}](N2) -- (M3);
    
    \draw [single={}](M3) -- (M4);
    \draw [double={}](M5) -- (M6);

    \end{tikzpicture}
    \captionsetup{width=1.0\linewidth}
  \captionof{figure}{
     the weighted forest from figure \ref{fig:3weighted}, supposing we were to remove a factor of 36.
         }
  \label{fig:3weighteddivided}
 \end{center}

\vfill\null
\columnbreak
 \begin{center}
 
        \tikzset{
        unmarked/.style args={}{
            draw,fill = black,circle,inner sep = 0pt, minimum size = .30cm,
        }
    }

    \tikzset{
        marked/.style args={}{
           draw, circle, fill=none, inner sep=0pt, minimum size=0.30cm, line width=0.04cm
        }
    }
    
    \tikzset{
        marked2/.style args={}{
           draw, double, circle, fill=none, inner sep=0pt, minimum size=0.30cm, line width=0.04cm
        }
    }

    \begin{tikzpicture}[scale=1]
        
        \node [marked={},label=below:{\small5}] (M1) at (-2.5,0) {};
        
        \node [unmarked={}, label=below:{\small6}] (M2) at (-1.5,0) {};
        \node [marked={},   label=below:{\small7}] (M3) at (-0.5,0) {};
        \node [unmarked={}, label=below:{\small8}] (M4) at (0.5,0) {};
        \node [unmarked={}, label=below:{\small9}] (M5) at (1.5,0) {};
        \node [marked={},   label=below:{\small10}] (M6) at (2.5,0) {};
        \node [unmarked={}, label=above:{\small3}] (N1) at (0.5,1.5) {};
        \node [marked={},   label=above:{\small4}] (N2) at (1.5,1.5) {};
        
        \node [marked={},label=above:{\small1}] (C1) at (-1.5,1.5) {};
        \node [marked={},label=above:{\small2}] (C2) at (-0.5,1.5) {};
    
    \draw (N1) -- (M2);
    \draw (N1) -- (M4);
    \draw (N2) -- (M3);
    
    \draw (M3) -- (M4);
    \draw (M5) -- (M6);

    \end{tikzpicture}
       \captionsetup{width=1.0\linewidth}
  \captionof{figure}{The trirooted forest produced by using the reverse algorithm.}
  \label{fig:trirooteddivided}
 \end{center}
 
 \vfill\null
\columnbreak
 \begin{center}
 
        \tikzset{
        unmarked/.style args={}{
            draw,fill = black,circle,inner sep = 0pt, minimum size = .30cm,
        }
    }

    \tikzset{
        marked/.style args={}{
           draw, circle, fill=none, inner sep=0pt, minimum size=0.30cm, line width=0.04cm
        }
    }
    
    \tikzset{
        marked2/.style args={}{
           draw, double, circle, fill=none, inner sep=0pt, minimum size=0.30cm, line width=0.04cm
        }
    }

    \begin{tikzpicture}[scale=1]
        
        \node [marked={},label=below:{\small5}] (M1) at (-2.5,0) {};
        
        \node [unmarked={}, label=below:{\small6}] (M2) at (-1.5,0) {};
        \node [marked={}, label=below:{\small7}] (M3) at (-0.5,0) {};
        \node [unmarked={},   label=below:{\small8}] (M4) at (0.5,0) {};
        \node [unmarked={}, label=below:{\small9}] (M5) at (1.5,0) {};
        \node [marked={},   label=below:{\small10}] (M6) at (2.5,0) {};
        \node [unmarked={}, label=above:{\small3}] (N1) at (0.5,1.5) {};
        \node [marked={},   label=above:{\small4}] (N2) at (1.5,1.5) {};
        
        \node [unmarked={},label=above:{\small0}] (C0) at (0,3) {};
        \node [marked={},label=above:{\small1}] (C1) at (-1.5,1.5) {};
        \node [marked={},label=above:{\small2}] (C2) at (-0.5,1.5) {};
    
    \draw [red](C1) -- (M1);
    \draw [red](C1) -- (C0);
    \draw [red](C2) -- (M3);
    \draw [red](C2) -- (M6);

    \draw [red](C0) -- (N2);
    \draw (N1) -- (M2);
    \draw (N1) -- (M4);
    \draw (N2) -- (M3);
    
    \draw (M3) -- (M4);
    \draw (M5) -- (M6);

    \end{tikzpicture}
      \captionsetup{width=1.0\linewidth}
  \captionof{figure}{The associated 3-weighted forest.}
    \label{fig:spantreedivided}
  \end{center}

\end{multicols}

\end{figure}

\begin{proposition}
The image of $\phi$ is pure.
\end{proposition}
\begin{proof}
Recall that from Corollary~\ref{cor:maxpassivespanning}, a monomial is of largest degree if and only if all normal marks are passive and the blueprint is spanning. Hence we need to show that given any trirooted forest $T$ that does not satisfy this property, we can come up with a new $3$-weighted forest $(F',w')$ that is divisible by $\phi_2(T) = (F,w)$. 


The only case when we have a normal mark on a vertex $v$ that is not passive, is when the mark is the sole mark of the component and is located at the smallest vertex of the component. In $(F,w)$ the component $C$ containing $v$ has no empowered edges. If $C$ only consists of $v$, pick any edge $e$ adjacent to $v$ in $\GG$ (which exists since we are assuming no coloops from Remark~\ref{rem:noloopcoloop}) and add it to $F$ with weight $1$. This component is obviously correctly weighted and the blueprint doesn't change, so we get the desired $3$-weighted forest. If $C$ has at least one edge, empower any arbitrary edge of this component. We still get a correctly weighted component, and the blueprint stays the same so we get a $3$-weighted forest we desire. 

Hence from now on we may assume that all normal marks of $T$ are passive. This means that only potential non-passive edges in the corresponding spanning tree $\phi_1^{-1}(T)$ are $01$ and $02$. Let us first start with the case when $01$ is an edge of $T$. From Remark~\ref{rem:noloopcoloop}, since $01$ is not a coloop, we can replace $01$ some other edge $e$ to get a spanning tree of $G$. If there is a path from $1$ to $2$ in the graph not using $01$ (which corresponds to $\P(\GG)$ containing the edge $12$), then we can set $e = 02$. The resulting trirooted forest is obtained from $T$ by adding an extra mark to $2$ and making it a doubly marked vertex. This new mark will empower some edge and keep the previous empowerments coming from other marks. In the original blueprint, it contained the edge $12$ and didn't contain $01$ or $02$ thanks to Lemma~\ref{lem:01blueprint}. The new blueprint consists of edges $01$ and $02$ so is acyclic. Hence this is the desired $3$-weighted forest.

The remaining case to consider is when there was no path from $1$ to $2$ in the graph not using $01$, and the components of $\GG$ containing a mark of type $1$ has no other marks in it. Moreover any regular mark is passive. In terms of the blueprint, thanks to Lemma~\ref{lem:01blueprint}, it is either empty or only consists of the edge $02$.

If there is a component $C$ with a sole mark of type $1$, there is an empowered edge $f$ coming from that mark in $C$. We wish to show that if there exists a vertex of type $\not = 1$ in $C$, then we can add an empowerment while keeping it a correctly weighted tree. If one of the endpoint of $f$ is not of type $1$, we simply increase $w_C(f)$ from $2$ to $3$. If not, look at the minimal path from the vertex of type $\not = 1$ to $e$ and empower the first edge along this path. We get a correctly weighted tree, and the blueprint gets an edge adjacent to $1$ added but is still acyclic. Hence this is the desired $3$-weighted forest.

Therefore we may assume that all components containing a mark of type $1$ only contains vertices of type $1$. Now assume that a component $C$ with type set $\subseteq \{0,2\}$ contains a vertex $v$ of type $1$. Unless in the case we have two marks $a,b$ in $C$ each of type $0$ and $2$ respectively, and $a$ lies in a path between $v$ and $b$, we simply add a mark to $v$ to get a correctly weighted tree such that its image under $\phi_2$ is a correctly weighted tree and divides $\phi_2(C)$. In terms of the blueprint, we are adding an edge adjacent to $1$ and is acyclic. If in that case, if there is a vertex of type $1$ in the path between $a$ and $b$, set that as $v$ instead. Otherwise, we put a mark on $v$ and at the same time move the mark from $a$ to the first vertex along the path from $a$ to $b$ that is of type $0$. If all vertices as we go from $a$ to $b$ are of type $2$ (aside from $a$), we move the mark of $a$ to $b$ and make $b$ a doubly-marked vertex. This gives a correctly marked tree such that its image under $\phi_2$ is a correctly weighted tree and is divisible by $\phi_2(C)$. In terms of the blueprint, we are adding an edge adjacent to $1$ and is acyclic. In any of those two cases, we get a $3$-weighted forest we desire.

Hence we may assume that all components containing a mark of type $1$ only consists of vertices of type $1$, and other components only contain vertices of type $0$ or $2$. So the edge $e$ that we can replace $01$ with will be an edge outside $B_0$ and adding it will merge a component with vertices of type $1$ and a component $C$ with vertices of type $0$ or $2$. Thanks to this condition compatibility between empowered objects on different sides is guaranteed. So the resulting merged component is again a correctly weighted component and we get the desired $3$-weighted forest.

Finally we have to consider the case when the spanning tree $B$ doesn't contain $01$ but contains $02$ as an active edge. Since $02$ is active it means there is no path from $2$ to $1$ that doesn't go through $02$ in $B$. From this point we pretty much replicate the proof for $01$ above: we first reduce the problem to the case when any component containing a mark of type $2$ can only contain vertices of that type. Then we further reduce to the case when any component not containing a mark of type $2$ can only contain vertices of type $0$ and $1$ (this part actually turns out to be much simpler that that for $01$ above. This is due to $2$ being bigger than both $0$ or $1$). So the edge $e$ that we can replace $02$ with to get another spanning tree of $G$, is going be one that connects a component consisting only of vertices of type $2$ to a component consisting only of vertices of type $\not = 2$. The merged component is correctly weighted and we get the desired $3$-weighted forest.

\end{proof}

For example, suppose we have the non-maximal spanning tree shown in figure \ref{fig:purityexa}, with its $3$-rooted forest and $3$-weighted forest in figures \ref{fig:purityexb} and \ref{fig:purityexc} respectively. Then, we can use the method used in the proof of purity above to construct a tree with higher passivity; since $9$ is an active normal mark, we replace it with a mark on the vertex $10$ as shown in figure \ref{fig:purityexd}. (For reference, figures \ref{fig:purityexe} and \ref{fig:purityexf} show the corresponding 3-weighted forest and spanning tree respectively.) $5$ is also an active normal mark, so we remove the mark on $5$ and add the edge $56$, shown in figure \ref{fig:purityexg}. (As before, see figures \ref{fig:purityexh} and \ref{fig:purityexi} for the corresponding 3-weighted forest and spanning tree respectively.) Finally, since $01$ is present, as there is a path from vertex $1$ to vertex $2$ which does not contain $01$, we replace $01$ with $02$. This can be seen in figure \ref{fig:purityexj}, with figure \ref{fig:purityexk} showing the final 3-weighted forest and \ref{fig:purityexl} the final spanning tree.

\begin{figure}[h]
\begin{multicols}{3}
 \begin{center}
 
        \tikzset{
        unmarked/.style args={}{
            draw,fill = black,circle,inner sep = 0pt, minimum size = .30cm,
        }
    }

    \tikzset{
        marked/.style args={}{
           draw, circle, fill=none, inner sep=0pt, minimum size=0.30cm, line width=0.04cm
        }
    }
    
    \tikzset{
        marked2/.style args={}{
           draw, double, circle, fill=none, inner sep=0pt, minimum size=0.30cm, line width=0.04cm
        }
    }

    \begin{tikzpicture}[scale=1]
        
        \node [marked={},label=below:{\small5}] (M1) at (-2.5,0) {};
        
        \node [marked={}, label=below:{\small6}] (M2) at (-1.5,0) {};
        \node [unmarked={}, label=below:{\small7}] (M3) at (-0.5,0) {};
        \node [marked={},   label=below:{\small8}] (M4) at (0.5,0) {};
        \node [marked={}, label=below:{\small9}] (M5) at (1.5,0) {};
        \node [unmarked={},   label=below:{\small10}] (M6) at (2.5,0) {};
        \node [unmarked={}, label=above:{\small3}] (N1) at (0.5,1.5) {};
        \node [unmarked={},   label=above:{\small4}] (N2) at (1.5,1.5) {};
        
        \node [unmarked={},label=above:{\small0}] (C0) at (0,3) {};
        \node [marked={},label=above:{\small1}] (C1) at (-1.5,1.5) {};
        \node [marked={},label=above:{\small2}] (C2) at (-0.5,1.5) {};
    
    \draw [red](C1) -- (M1);
    \draw [red](C1) -- (M2);
    \draw [red](C2) -- (M4);
    \draw [red](C2) -- (M5);

    \draw [red](C0) -- (C1);
    \draw (N1) -- (M2);
    \draw (N1) -- (M4);
    \draw (N2) -- (M3);
    
    \draw (M3) -- (M4);
    \draw (M5) -- (M6);

    \end{tikzpicture}
    \captionsetup{width=1.0\linewidth}
  \captionof{figure}{A non-maximal spanning tree}
  \label{fig:purityexa}
 \end{center}

\vfill\null
\columnbreak
 \begin{center}
 
        \tikzset{
        unmarked/.style args={}{
            draw,fill = black,circle,inner sep = 0pt, minimum size = .30cm,
        }
    }

    \tikzset{
        marked/.style args={}{
           draw, circle, fill=none, inner sep=0pt, minimum size=0.30cm, line width=0.04cm
        }
    }
    
    \tikzset{
        marked2/.style args={}{
           draw, double, circle, fill=none, inner sep=0pt, minimum size=0.30cm, line width=0.04cm
        }
    }

    \begin{tikzpicture}[scale=1]
        
        \node [marked={},label=below:{\small5}] (M1) at (-2.5,0) {};
        
        \node [marked={}, label=below:{\small6}] (M2) at (-1.5,0) {};
        \node [unmarked={}, label=below:{\small7}] (M3) at (-0.5,0) {};
        \node [marked={},   label=below:{\small8}] (M4) at (0.5,0) {};
        \node [marked={}, label=below:{\small9}] (M5) at (1.5,0) {};
        \node [unmarked={},   label=below:{\small10}] (M6) at (2.5,0) {};
        \node [unmarked={}, label=above:{\small3}] (N1) at (0.5,1.5) {};
        \node [unmarked={},   label=above:{\small4}] (N2) at (1.5,1.5) {};
        
        \node [marked={},label=above:{\small1}] (C1) at (-1.5,1.5) {};
        \node [marked={},label=above:{\small2}] (C2) at (-0.5,1.5) {};

    \draw (N1) -- (M2);
    \draw (N1) -- (M4);
    \draw (N2) -- (M3);
    
    \draw (M3) -- (M4);
    \draw (M5) -- (M6);

    \end{tikzpicture}
       \captionsetup{width=1.0\linewidth}
  \captionof{figure}{The corresponding 3-rooted forest}
  \label{fig:purityexb}
 \end{center}
 
 \vfill\null
\columnbreak
 \begin{center}
 
        \tikzset{
        unmarked/.style args={}{
            draw,fill = black,circle,inner sep = 0pt, minimum size = .30cm,
        }
    }

    \tikzset{
        marked/.style args={}{
           draw, circle, fill=none, inner sep=0pt, minimum size=0.30cm, line width=0.04cm
        }
    }
    
    \tikzset{
        marked2/.style args={}{
           draw, double, circle, fill=none, inner sep=0pt, minimum size=0.30cm, line width=0.04cm
        }
    }
    
    \tikzset{
        single/.style args={}{
            draw,line width=0.5mm,black
        }
    }
    \tikzset{
        double/.style args={}{
            preaction={draw,line width=1.25mm,black},draw,line width=.25mm,white
        }
    }
    \tikzset{
        triple/.style args={}{
            preaction={preaction={draw,line width=2mm,black},draw,line width=1mm,white},draw,line 
            width=.5mm,black
        }
    }  
    \tikzset{
        quadruple/.style args={}{
            preaction={preaction={preaction={draw,line width=2.75mm,black},draw,line width=1.75mm,white},draw,line width=1.25mm,black},draw,line width=.25mm,white
        }
    }

    \begin{tikzpicture}[scale=1]
        
        \node [unmarked={},label=below:{\small5}] (M1) at (-2.5,0) {};
        
        \node [unmarked={}, label=below:{\small6}] (M2) at (-1.5,0) {};
        \node [unmarked={}, label=below:{\small7}] (M3) at (-0.5,0) {};
        \node [unmarked={},   label=below:{\small8}] (M4) at (0.5,0) {};
        \node [unmarked={}, label=below:{\small9}] (M5) at (1.5,0) {};
        \node [unmarked={},   label=below:{\small10}] (M6) at (2.5,0) {};
        \node [unmarked={}, label=above:{\small3}] (N1) at (0.5,1.5) {};
        \node [unmarked={},   label=above:{\small4}] (N2) at (1.5,1.5) {};
        
        \node [unmarked={},label=above:{\small1}] (C1) at (-1.5,1.5) {};
        \node [unmarked={},label=above:{\small2}] (C2) at (-0.5,1.5) {};

    \draw [double={}](N1) -- (M2);
    \draw [double={}](N1) -- (M4);
    \draw [single={}](N2) -- (M3);
    
    \draw [single={}](M3) -- (M4);
    \draw [single={}](M5) -- (M6);

    \end{tikzpicture}
      \captionsetup{width=1.0\linewidth}
  \captionof{figure}{The corresponding 3-weighted forest}
    \label{fig:purityexc}
    \end{center}

\end{multicols}
\end{figure}

\begin{figure}[h]
\begin{multicols}{3}
 \begin{center}

  \vspace{5mm}
        \tikzset{
        unmarked/.style args={}{
            draw,fill = black,circle,inner sep = 0pt, minimum size = .30cm,
        }
    }

    \tikzset{
        marked/.style args={}{
           draw, circle, fill=none, inner sep=0pt, minimum size=0.30cm, line width=0.04cm
        }
    }
    
    \tikzset{
        marked2/.style args={}{
           draw, double, circle, fill=none, inner sep=0pt, minimum size=0.30cm, line width=0.04cm
        }
    }

    \begin{tikzpicture}[scale=1]
        
        \node [marked={},label=below:{\small5}] (M1) at (-2.5,0) {};
        
        \node [marked={}, label=below:{\small6}] (M2) at (-1.5,0) {};
        \node [unmarked={}, label=below:{\small7}] (M3) at (-0.5,0) {};
        \node [marked={},   label=below:{\small8}] (M4) at (0.5,0) {};
        \node [unmarked={}, label=below:{\small9}] (M5) at (1.5,0) {};
        \node [marked={},   label=below:{\small10}] (M6) at (2.5,0) {};
        \node [unmarked={}, label=above:{\small3}] (N1) at (0.5,1.5) {};
        \node [unmarked={},   label=above:{\small4}] (N2) at (1.5,1.5) {};
        
        \node [marked={},label=above:{\small1}] (C1) at (-1.5,1.5) {};
        \node [marked={},label=above:{\small2}] (C2) at (-0.5,1.5) {};

    \draw (N1) -- (M2);
    \draw (N1) -- (M4);
    \draw (N2) -- (M3);
    
    \draw (M3) -- (M4);
    \draw (M5) -- (M6);

    \end{tikzpicture}
    \captionsetup{width=1.0\linewidth}
  \captionof{figure}{Moving the mark from vertex $9$ to vertex $10$}
  \label{fig:purityexd}
  
 \end{center}

\vfill\null
\columnbreak
 \begin{center}

  \vspace{20mm}
        \tikzset{
        unmarked/.style args={}{
            draw,fill = black,circle,inner sep = 0pt, minimum size = .30cm,
        }
    }

    \tikzset{
        marked/.style args={}{
           draw, circle, fill=none, inner sep=0pt, minimum size=0.30cm, line width=0.04cm
        }
    }
    
    \tikzset{
        marked2/.style args={}{
           draw, double, circle, fill=none, inner sep=0pt, minimum size=0.30cm, line width=0.04cm
        }
    }
\tikzset{
        single/.style args={}{
            draw,line width=0.5mm,black
        }
    }
    \tikzset{
        double/.style args={}{
            preaction={draw,line width=1.25mm,black},draw,line width=.25mm,white
        }
    }
    \tikzset{
        triple/.style args={}{
            preaction={preaction={draw,line width=2mm,black},draw,line width=1mm,white},draw,line 
            width=.5mm,black
        }
    }  
    \tikzset{
        quadruple/.style args={}{
            preaction={preaction={preaction={draw,line width=2.75mm,black},draw,line width=1.75mm,white},draw,line width=1.25mm,black},draw,line width=.25mm,white
        }
    }
    \begin{tikzpicture}[scale=1]
        
        \node [unmarked={},label=below:{\small5}] (M1) at (-2.5,0) {};
        
        \node [unmarked={}, label=below:{\small6}] (M2) at (-1.5,0) {};
        \node [unmarked={}, label=below:{\small7}] (M3) at (-0.5,0) {};
        \node [unmarked={},   label=below:{\small8}] (M4) at (0.5,0) {};
        \node [unmarked={}, label=below:{\small9}] (M5) at (1.5,0) {};
        \node [unmarked={},   label=below:{\small10}] (M6) at (2.5,0) {};
        \node [unmarked={}, label=above:{\small3}] (N1) at (0.5,1.5) {};
        \node [unmarked={},   label=above:{\small4}] (N2) at (1.5,1.5) {};
        
        \node [unmarked={},label=above:{\small1}] (C1) at (-1.5,1.5) {};
        \node [unmarked={},label=above:{\small2}] (C2) at (-0.5,1.5) {};

    \draw [double={}](N1) -- (M2);
    \draw [double={}](N1) -- (M4);
    \draw [single={}](N2) -- (M3);
    
    \draw [single={}](M3) -- (M4);
    \draw [double={}](M5) -- (M6);

    \end{tikzpicture}
    \captionsetup{width=1.0\linewidth}
  \captionof{figure}{The corresponding $3$-weighted forest. The same forest as before, but now with an extra weight on $9A$.}
  \label{fig:purityexe}
  
 \end{center}
 
 \vfill\null
\columnbreak
 \begin{center}

  \vspace{5mm}
        \tikzset{
        unmarked/.style args={}{
            draw,fill = black,circle,inner sep = 0pt, minimum size = .30cm,
        }
    }

    \tikzset{
        marked/.style args={}{
           draw, circle, fill=none, inner sep=0pt, minimum size=0.30cm, line width=0.04cm
        }
    }
    
    \tikzset{
        marked2/.style args={}{
           draw, double, circle, fill=none, inner sep=0pt, minimum size=0.30cm, line width=0.04cm
        }
    }

    \begin{tikzpicture}[scale=1]
        
        \node [marked={},label=below:{\small5}] (M1) at (-2.5,0) {};
        
        \node [marked={}, label=below:{\small6}] (M2) at (-1.5,0) {};
        \node [unmarked={}, label=below:{\small7}] (M3) at (-0.5,0) {};
        \node [marked={},   label=below:{\small8}] (M4) at (0.5,0) {};
        \node [unmarked={}, label=below:{\small9}] (M5) at (1.5,0) {};
        \node [marked={},   label=below:{\small10}] (M6) at (2.5,0) {};
        \node [unmarked={}, label=above:{\small3}] (N1) at (0.5,1.5) {};
        \node [unmarked={},   label=above:{\small4}] (N2) at (1.5,1.5) {};
        
        \node [unmarked={},label=above:{\small0}] (C0) at (0,3) {};
        \node [marked={},label=above:{\small1}] (C1) at (-1.5,1.5) {};
        \node [marked={},label=above:{\small2}] (C2) at (-0.5,1.5) {};
    
    \draw [red](C1) -- (M1);
    \draw [red](C1) -- (M2);
    \draw [red](C2) -- (M4);
    \draw [red](C2) -- (M6);

    \draw [red](C0) -- (C1);
    \draw (N1) -- (M2);
    \draw (N1) -- (M4);
    \draw (N2) -- (M3);
    
    \draw (M3) -- (M4);
    \draw (M5) -- (M6);

    \end{tikzpicture}
    \captionsetup{width=1.0\linewidth}
  \captionof{figure}{The corresponding spanning tree}
  \label{fig:purityexf}
  \end{center}

\end{multicols}
\end{figure}

\begin{figure}[h]
\begin{multicols}{3}
 \begin{center}

  \vspace{20mm}
        \tikzset{
        unmarked/.style args={}{
            draw,fill = black,circle,inner sep = 0pt, minimum size = .30cm,
        }
    }

    \tikzset{
        marked/.style args={}{
           draw, circle, fill=none, inner sep=0pt, minimum size=0.30cm, line width=0.04cm
        }
    }
    
    \tikzset{
        marked2/.style args={}{
           draw, double, circle, fill=none, inner sep=0pt, minimum size=0.30cm, line width=0.04cm
        }
    }

    \begin{tikzpicture}[scale=1]
        
        \node [unmarked={},label=below:{\small5}] (M1) at (-2.5,0) {};
        
        \node [marked={}, label=below:{\small6}] (M2) at (-1.5,0) {};
        \node [unmarked={}, label=below:{\small7}] (M3) at (-0.5,0) {};
        \node [marked={},   label=below:{\small8}] (M4) at (0.5,0) {};
        \node [unmarked={}, label=below:{\small9}] (M5) at (1.5,0) {};
        \node [marked={},   label=below:{\small10}] (M6) at (2.5,0) {};
        \node [unmarked={}, label=above:{\small3}] (N1) at (0.5,1.5) {};
        \node [unmarked={},   label=above:{\small4}] (N2) at (1.5,1.5) {};
        
        \node [marked={},label=above:{\small1}] (C1) at (-1.5,1.5) {};
        \node [marked={},label=above:{\small2}] (C2) at (-0.5,1.5) {};

    \draw (N1) -- (M2);
    \draw (N1) -- (M4);
    \draw (N2) -- (M3);
    
    \draw (M1) -- (M2);
    \draw (M3) -- (M4);
    \draw (M5) -- (M6);

    \end{tikzpicture}
    \captionsetup{width=1.0\linewidth}
  \captionof{figure}{Replacing the mark on vertex $5$ with the edge $56$}
  \label{fig:purityexg}
 \end{center}

\vfill\null
\columnbreak
 \begin{center}

  \vspace{20mm}
        \tikzset{
        unmarked/.style args={}{
            draw,fill = black,circle,inner sep = 0pt, minimum size = .30cm,
        }
    }

    \tikzset{
        marked/.style args={}{
           draw, circle, fill=none, inner sep=0pt, minimum size=0.30cm, line width=0.04cm
        }
    }
    
    \tikzset{
        marked2/.style args={}{
           draw, double, circle, fill=none, inner sep=0pt, minimum size=0.30cm, line width=0.04cm
        }
    }
\tikzset{
        single/.style args={}{
            draw,line width=0.5mm,black
        }
    }
    \tikzset{
        double/.style args={}{
            preaction={draw,line width=1.25mm,black},draw,line width=.25mm,white
        }
    }
    \tikzset{
        triple/.style args={}{
            preaction={preaction={draw,line width=2mm,black},draw,line width=1mm,white},draw,line 
            width=.5mm,black
        }
    }  
    \tikzset{
        quadruple/.style args={}{
            preaction={preaction={preaction={draw,line width=2.75mm,black},draw,line width=1.75mm,white},draw,line width=1.25mm,black},draw,line width=.25mm,white
        }
    }
    \begin{tikzpicture}[scale=1]
        
        \node [unmarked={},label=below:{\small5}] (M1) at (-2.5,0) {};
        
        \node [unmarked={}, label=below:{\small6}] (M2) at (-1.5,0) {};
        \node [unmarked={}, label=below:{\small7}] (M3) at (-0.5,0) {};
        \node [unmarked={},   label=below:{\small8}] (M4) at (0.5,0) {};
        \node [unmarked={}, label=below:{\small9}] (M5) at (1.5,0) {};
        \node [unmarked={},   label=below:{\small10}] (M6) at (2.5,0) {};
        \node [unmarked={}, label=above:{\small3}] (N1) at (0.5,1.5) {};
        \node [unmarked={},   label=above:{\small4}] (N2) at (1.5,1.5) {};
        
        \node [unmarked={},label=above:{\small1}] (C1) at (-1.5,1.5) {};
        \node [unmarked={},label=above:{\small2}] (C2) at (-0.5,1.5) {};

    \draw [single={}](M1) -- (M2);
    \draw [double={}](N1) -- (M2);
    \draw [double={}](N1) -- (M4);
    \draw [single={}](N2) -- (M3);
    
    \draw [single={}](M3) -- (M4);
    \draw [double={}](M5) -- (M6);

    \end{tikzpicture}
    \captionsetup{width=1.0\linewidth}
  \captionof{figure}{The corresponding $3$-weighted forest. The edge $56$ has been added.}
  \label{fig:purityexh}

 \end{center}
 
 \vfill\null
\columnbreak
 \begin{center}
 
  \vspace{3.7mm}
        \tikzset{
        unmarked/.style args={}{
            draw,fill = black,circle,inner sep = 0pt, minimum size = .30cm,
        }
    }

    \tikzset{
        marked/.style args={}{
           draw, circle, fill=none, inner sep=0pt, minimum size=0.30cm, line width=0.04cm
        }
    }
    
    \tikzset{
        marked2/.style args={}{
           draw, double, circle, fill=none, inner sep=0pt, minimum size=0.30cm, line width=0.04cm
        }
    }

    \begin{tikzpicture}[scale=1]
        
        \node [unmarked={},label=below:{\small5}] (M1) at (-2.5,0) {};
        
        \node [marked={}, label=below:{\small6}] (M2) at (-1.5,0) {};
        \node [unmarked={}, label=below:{\small7}] (M3) at (-0.5,0) {};
        \node [marked={},   label=below:{\small8}] (M4) at (0.5,0) {};
        \node [unmarked={}, label=below:{\small9}] (M5) at (1.5,0) {};
        \node [marked={},   label=below:{\small10}] (M6) at (2.5,0) {};
        \node [unmarked={}, label=above:{\small3}] (N1) at (0.5,1.5) {};
        \node [unmarked={},   label=above:{\small4}] (N2) at (1.5,1.5) {};
        
        \node [unmarked={},label=above:{\small0}] (C0) at (0,3) {};
        \node [marked={},label=above:{\small1}] (C1) at (-1.5,1.5) {};
        \node [marked={},label=above:{\small2}] (C2) at (-0.5,1.5) {};
    
    \draw [red](C1) -- (M2);
    \draw [red](C2) -- (M4);
    \draw [red](C2) -- (M6);

    \draw [red](C0) -- (C1);
    \draw (N1) -- (M2);
    \draw (N1) -- (M4);
    \draw (N2) -- (M3);
    
    \draw (M1) -- (M2);
    \draw (M3) -- (M4);
    \draw (M5) -- (M6);

    \end{tikzpicture}
    \captionsetup{width=1.0\linewidth}
  \captionof{figure}{The corresponding spanning tree}
  \label{fig:purityexi}

\end{center}
\end{multicols}
\end{figure}

\begin{figure}[h]
\begin{multicols}{3}
 \begin{center}

  \vspace{20mm}
        \tikzset{
        unmarked/.style args={}{
            draw,fill = black,circle,inner sep = 0pt, minimum size = .30cm,
        }
    }

    \tikzset{
        marked/.style args={}{
           draw, circle, fill=none, inner sep=0pt, minimum size=0.30cm, line width=0.04cm
        }
    }
    
    \tikzset{
        marked2/.style args={}{
           draw, double, circle, fill=none, inner sep=0pt, minimum size=0.30cm, line width=0.04cm
        }
    }

    \begin{tikzpicture}[scale=1]
        
        \node [unmarked={},label=below:{\small5}] (M1) at (-2.5,0) {};
        
        \node [marked={}, label=below:{\small6}] (M2) at (-1.5,0) {};
        \node [unmarked={}, label=below:{\small7}] (M3) at (-0.5,0) {};
        \node [marked2={},   label=below:{\small8}] (M4) at (0.5,0) {};
        \node [unmarked={}, label=below:{\small9}] (M5) at (1.5,0) {};
        \node [marked={},   label=below:{\small10}] (M6) at (2.5,0) {};
        \node [unmarked={}, label=above:{\small3}] (N1) at (0.5,1.5) {};
        \node [unmarked={},   label=above:{\small4}] (N2) at (1.5,1.5) {};
        
        \node [marked={},label=above:{\small1}] (C1) at (-1.5,1.5) {};
        \node [marked={},label=above:{\small2}] (C2) at (-0.5,1.5) {};

    \draw (N1) -- (M2);
    \draw (N1) -- (M4);
    \draw (N2) -- (M3);
    
    \draw (M1) -- (M2);
    \draw (M3) -- (M4);
    \draw (M5) -- (M6);

    \end{tikzpicture}
    \captionsetup{width=1.0\linewidth}
  \captionof{figure}{Replacing the edge $01$ with the edge $02$}
  \label{fig:purityexj}
 \end{center}

\vfill\null
\columnbreak
 \begin{center}

  \vspace{20mm}
        \tikzset{
        unmarked/.style args={}{
            draw,fill = black,circle,inner sep = 0pt, minimum size = .30cm,
        }
    }

    \tikzset{
        marked/.style args={}{
           draw, circle, fill=none, inner sep=0pt, minimum size=0.30cm, line width=0.04cm
        }
    }
    
    \tikzset{
        marked2/.style args={}{
           draw, double, circle, fill=none, inner sep=0pt, minimum size=0.30cm, line width=0.04cm
        }
    }
\tikzset{
        single/.style args={}{
            draw,line width=0.5mm,black
        }
    }
    \tikzset{
        double/.style args={}{
            preaction={draw,line width=1.25mm,black},draw,line width=.25mm,white
        }
    }
    \tikzset{
        triple/.style args={}{
            preaction={preaction={draw,line width=2mm,black},draw,line width=1mm,white},draw,line 
            width=.5mm,black
        }
    }  
    \tikzset{
        quadruple/.style args={}{
            preaction={preaction={preaction={draw,line width=2.75mm,black},draw,line width=1.75mm,white},draw,line width=1.25mm,black},draw,line width=.25mm,white
        }
    }
    \begin{tikzpicture}[scale=1]
        
        \node [unmarked={},label=below:{\small5}] (M1) at (-2.5,0) {};
        
        \node [unmarked={}, label=below:{\small6}] (M2) at (-1.5,0) {};
        \node [unmarked={}, label=below:{\small7}] (M3) at (-0.5,0) {};
        \node [unmarked={},   label=below:{\small8}] (M4) at (0.5,0) {};
        \node [unmarked={}, label=below:{\small9}] (M5) at (1.5,0) {};
        \node [unmarked={},   label=below:{\small10}] (M6) at (2.5,0) {};
        \node [unmarked={}, label=above:{\small3}] (N1) at (0.5,1.5) {};
        \node [unmarked={},   label=above:{\small4}] (N2) at (1.5,1.5) {};
        
        \node [unmarked={},label=above:{\small1}] (C1) at (-1.5,1.5) {};
        \node [unmarked={},label=above:{\small2}] (C2) at (-0.5,1.5) {};

    \draw [single={}](M1) -- (M2);
    \draw [double={}](N1) -- (M2);
    \draw [triple={}](N1) -- (M4);
    \draw [single={}](N2) -- (M3);
    
    \draw [single={}](M3) -- (M4);
    \draw [double={}](M5) -- (M6);

    \end{tikzpicture}
    \captionsetup{width=1.0\linewidth}
  \captionof{figure}{The corresponding $3$-weighted forest. $38$ now has weight 3.}
  \label{fig:purityexk}
 \end{center}
 
 \vfill\null
\columnbreak
 \begin{center}
 
  \vspace{3.7mm}
        \tikzset{
        unmarked/.style args={}{
            draw,fill = black,circle,inner sep = 0pt, minimum size = .30cm,
        }
    }

    \tikzset{
        marked/.style args={}{
           draw, circle, fill=none, inner sep=0pt, minimum size=0.30cm, line width=0.04cm
        }
    }
    
    \tikzset{
        marked2/.style args={}{
           draw, double, circle, fill=none, inner sep=0pt, minimum size=0.30cm, line width=0.04cm
        }
    }

    \begin{tikzpicture}[scale=1]
        
        \node [unmarked={},label=below:{\small5}] (M1) at (-2.5,0) {};
        
        \node [marked={},   label=below:{\small6}] (M2) at (-1.5,0) {};
        \node [unmarked={}, label=below:{\small7}] (M3) at (-0.5,0) {};
        \node [marked2={},  label=below:{\small8}] (M4) at (0.5,0) {};
        \node [unmarked={}, label=below:{\small9}] (M5) at (1.5,0) {};
        \node [marked={},   label=below:{\small10}] (M6) at (2.5,0) {};
        \node [unmarked={}, label=above:{\small3}] (N1) at (0.5,1.5) {};
        \node [unmarked={}, label=above:{\small4}] (N2) at (1.5,1.5) {};
        
        \node [unmarked={},label=above:{\small0}] (C0) at (0,3) {};
        \node [marked={},label=above:{\small1}] (C1) at (-1.5,1.5) {};
        \node [marked={},label=above:{\small2}] (C2) at (-0.5,1.5) {};
    
    \draw [red](C1) -- (M2);
    \draw [red](C2) -- (M4);
    \draw [red](C2) -- (M6);

    \draw [red](C0) -- (C2);
    \draw (N1) -- (M2);
    \draw (N1) -- (M4);
    \draw (N2) -- (M3);
    
    \draw (M1) -- (M2);
    \draw (M3) -- (M4);
    \draw (M5) -- (M6);

    \end{tikzpicture}
    \captionsetup{width=1.0\linewidth}
  \captionof{figure}{The corresponding spanning tree}
  \label{fig:purityexl}
\end{center}
\end{multicols}
\end{figure}

\begin{proposition}
The map $\phi$ sends passivity (of a spanning tree) to degree (of a monomial).
\end{proposition}
\begin{proof}
From Lemma~\ref{lem:edgepassive}, the passivity of a spanning tree gets sent to the number of edges in the forest plus the number of passive marks under the map $\phi_1$. The underlying forest doesn't change under the map $\phi_2$, and each empowerment comes from a passive mark from Remark~\ref{rem:markempower}. Hence the degree of the resulting monomial is same as the passivity of the spanning tree we started out with.
\end{proof}

For example take a look at the spanning tree $B$ back in Figure~\ref{fig:spantreedivided}. This spanning tree has $04,27,2A,36,38,47,78,9A$ as passive edges. In the trirooted forest we get under the map $\phi_1$, the edges $36,38,47,78,9A$ are still edges in the forest and other edges become passive marks: marks at $4,7,A$. Notice that $1,2$ are special marks and $5$ is the mark located on the smallest (actually the unique) vertex of the component. In the $3$-weighted forest we get under the map $\phi_2$, the passive marks at $4,7,A$ each correspond to empowerment of edges $47,47,9A$. Hence we started with $8$ passivity (of a spanning tree) and ended up with a monomial $x_{36}x_{38}x_{47}^3 x_{78} x_{9A}^2$ which has degree $8$ as well.

Combining the above three propositions we get the desired result of this paper:

\begin{corollary}
Let $G$ be a triconed graph. Its $h$-vector is a pure $O$-sequence.
\end{corollary}


\section{Worked example and other discussions}

\subsection{Worked example}
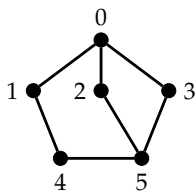
\begin{figure}[H]
  \begin{center}
 
  \scalebox{.75}{        \begin{tikzpicture}[scale=.6]
        \node [draw,fill,circle,inner sep = 0pt, minimum size = .25cm,label=above:{\large0}] (0) at (0,1.5) {};
        \node [draw,fill,circle,inner sep = 0pt, minimum size = .25cm,label=left:{\large1}] (1) at (-2,0) {};
        \node [draw,fill,circle,inner sep = 0pt, minimum size = .25cm,label=left:{\large2}] (2) at (0,0) {};
        \node [draw,fill,circle,inner sep = 0pt, minimum size = .25cm,label=right:{\large3}] (3) at (2,0) {};
        \node [draw,fill,circle,inner sep = 0pt, minimum size = .25cm,label=below:{\large4}] (4) at (-1.2,-2) {};
        \node [draw,fill,circle,inner sep = 0pt, minimum size = .25cm,label=below:{\large5}] (5) at (1.2,-2) {};
        
    \draw [line width = .5mm, black](0) -- (1);  
    \draw [line width = .5mm, black](0) -- (2);
    \draw [line width = .5mm, black](0) -- (3);  
    \draw [line width = .5mm, black](1) -- (4);  
    \draw [line width = .5mm, black](2) -- (5);  
    \draw [line width = .5mm, black](4) -- (5);  
    \draw [line width = .5mm, black](3) -- (5);

    \end{tikzpicture}}
     \captionsetup{width=1.0\linewidth}
  \captionof{figure}{An example graph $G$.}
  \label{fig:rungraph}
  \end{center}
 \end{figure}

Look at the triconed graph $G$ of Figure~\ref{fig:rungraph}. All the spanning trees of this graph are drawn in Figure~\ref{fig:runtrees}. The trirooted forests of the spanning trees we get by the map $\phi_1$ are drawn in Figure~\ref{fig:runtriiroot}. The $3$-weighted forests we get from the trirooted forests by the map $\phi_2$ are drawn in Figure~\ref{fig:runmulticomp}. In each figure, the poset structure is given by divisibility of the corresponding monomials. 

 \def\xSep{60}
 \def\ySep{120}
 
 \begin{figure}[H]
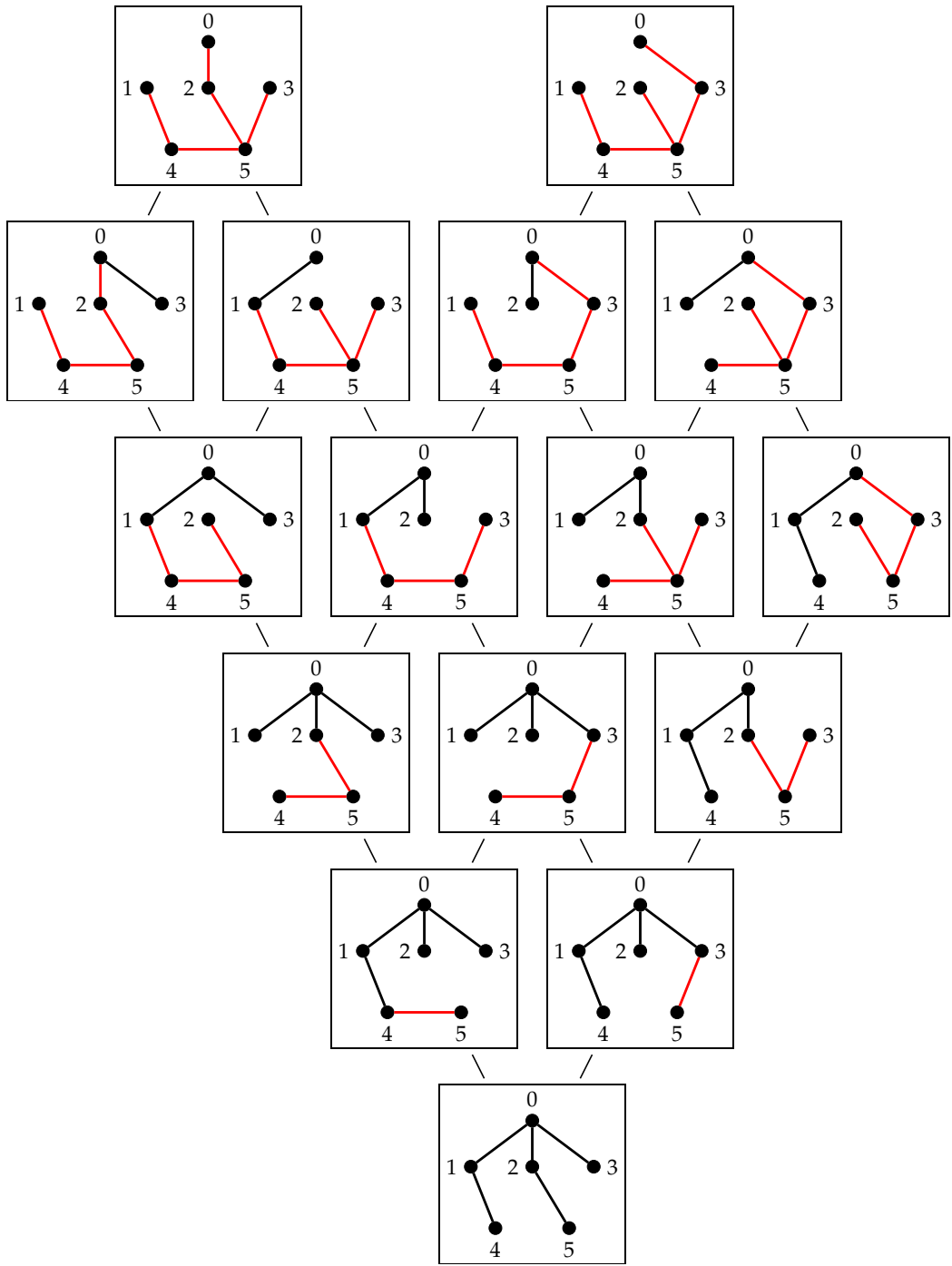

  \begin{center}
 
  \scalebox{.75}{        \begin{tikzpicture}[scale=.6]
        
        \node (N1) {\frame{\input{figures/posets/spanning trees/st1}}};
        
        \node[xshift=-1*\xSep, yshift=1*\ySep] (x) {\frame{\input{figures/posets/spanning trees/stx}}};
        \node[xshift=1*\xSep, yshift=1*\ySep] (y) {\frame{\input{figures/posets/spanning trees/sty}}};
        
        \node[xshift=-2*\xSep, yshift=2*\ySep] (x2) {\frame{\input{figures/posets/spanning trees/stx2}}};
        \node[xshift=0, yshift=2*\ySep] (xy) {\frame{\input{figures/posets/spanning trees/stxy}}};
        \node[xshift=2*\xSep, yshift=2*\ySep] (y2) {\frame{\input{figures/posets/spanning trees/sty2}}};
        
        \node[xshift=-3*\xSep, yshift=3*\ySep] (x3) {\frame{\input{figures/posets/spanning trees/stx3}}};
        \node[xshift=-1*\xSep, yshift=3*\ySep] (x2y) {\frame{\input{figures/posets/spanning trees/stx2y}}};
        \node[xshift=1*\xSep, yshift=3*\ySep] (xy2) {\frame{\input{figures/posets/spanning trees/stxy2}}};
        \node[xshift=3*\xSep, yshift=3*\ySep] (y3) {\frame{\input{figures/posets/spanning trees/sty3}}};
        
        \node[xshift=-4*\xSep, yshift=4*\ySep] (x4) {\frame{\input{figures/posets/spanning trees/stx4}}};
        \node[xshift=-2*\xSep, yshift=4*\ySep] (x3y) {\frame{\input{figures/posets/spanning trees/stx3y}}};
        \node[xshift=0, yshift=4*\ySep] (x2y2) {\frame{\input{figures/posets/spanning trees/stx2y2}}};
        \node[xshift=2*\xSep, yshift=4*\ySep] (xy3) {\frame{\input{figures/posets/spanning trees/stxy3}}};
        
        \node[xshift=-3*\xSep, yshift=5*\ySep] (x4y) {\frame{\input{figures/posets/spanning trees/stx4y}}};
        \node[xshift=1*\xSep, yshift=5*\ySep] (x2y3) {\frame{\input{figures/posets/spanning trees/stx2y3}}};
        
        \draw[thick] (N1) -- (x);
        \draw[thick] (N1) -- (y);
        
        \draw[thick] (x) -- (x2);
        \draw[thick] (x) -- (xy);
        \draw[thick] (y) -- (xy);
        \draw[thick] (y) -- (y2);
        
        \draw[thick] (x2) -- (x3);
        \draw[thick] (x2) -- (x2y);
        \draw[thick] (xy) -- (x2y);
        \draw[thick] (xy) -- (xy2);
        \draw[thick] (y2) -- (xy2);
        \draw[thick] (y2) -- (y3);
        
        \draw[thick] (x3) -- (x4);
        \draw[thick] (x3) -- (x3y);
        \draw[thick] (x2y) -- (x3y);
        \draw[thick] (x2y) -- (x2y2);
        \draw[thick] (xy2) -- (x2y2);
        \draw[thick] (xy2) -- (xy3);
        \draw[thick] (y3) -- (xy3);
        
        \draw[thick] (x4) -- (x4y);
        \draw[thick] (x3y) -- (x4y);
        \draw[thick] (x2y2) -- (x2y3);
        \draw[thick] (xy3) -- (x2y3);
    \end{tikzpicture}}
     \captionsetup{width=1.0\linewidth}
  \captionof{figure}{Spanning trees of $\G$. Passive edges are colored red.}
  \label{fig:runtrees}
  \end{center}
 \end{figure}


     \tikzset{
         unmarked/.style args={}{
             draw,fill = black,circle,inner sep = 0pt, minimum size = .30cm,
         }
     }

     \tikzset{
         marked/.style args={}{
          draw, circle, fill=none, inner sep=0pt, minimum size=0.30cm, line width=0.04cm
         }
     }
    
     \tikzset{
         marked2/.style args={}{
          draw, double, circle, fill=none, inner sep=0pt, minimum size=0.30cm, line width=0.04cm
         }
     }

 \def\xSep{50}
 \def\ySep{80}

 \begin{figure}[H]
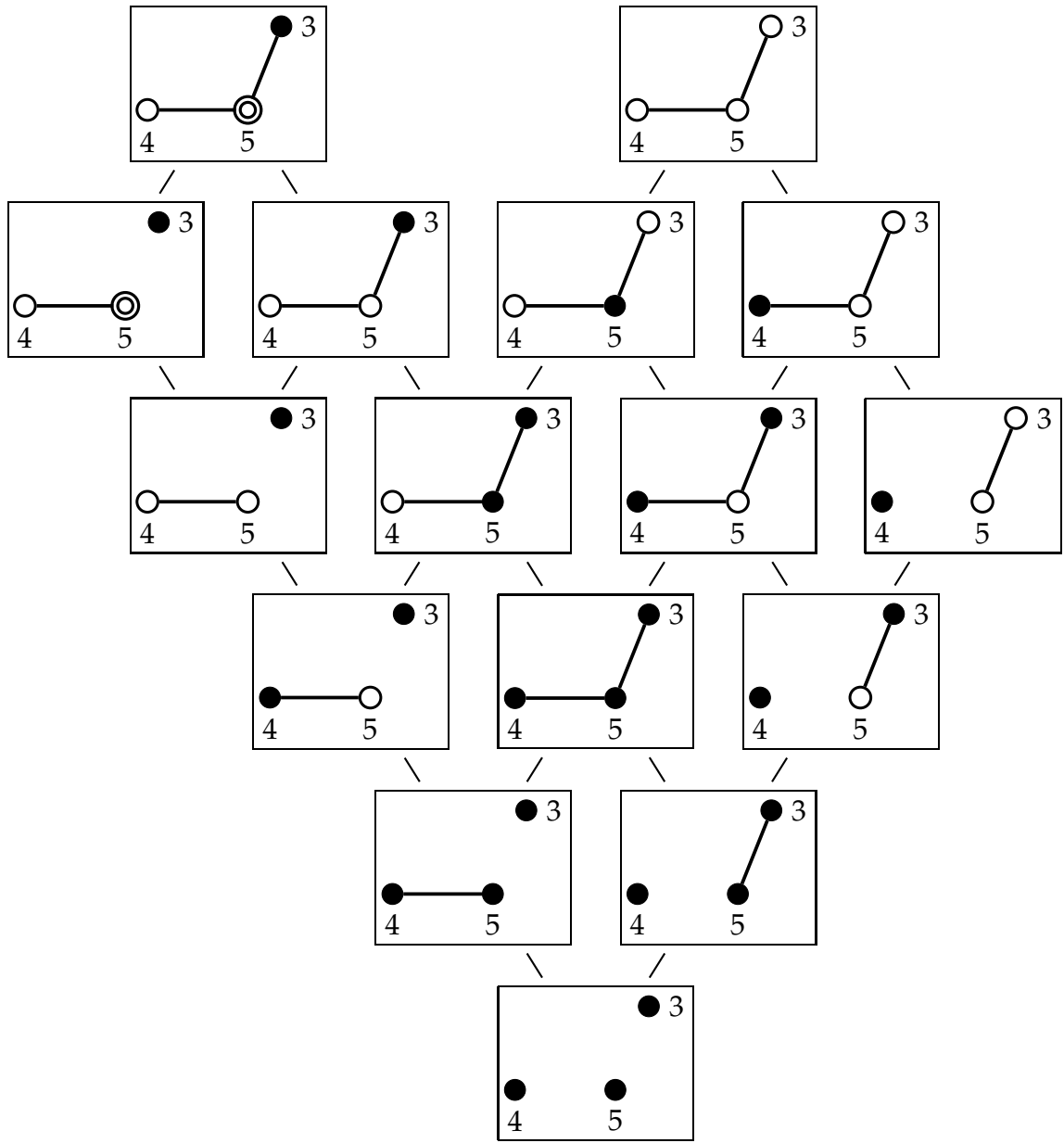

  \begin{center}
 
  \scalebox{1}{        \begin{tikzpicture}[scale=.6]
        
        \node (N1) {\frame{\input{figures/posets/vertex rooted/vr1}}};
        
        \node[xshift=-1*\xSep, yshift=1*\ySep] (x) {\frame{\input{figures/posets/vertex rooted/vrx}}};
        \node[xshift=1*\xSep, yshift=1*\ySep] (y) {\frame{\input{figures/posets/vertex rooted/vry}}};
        
        \node[xshift=-2*\xSep, yshift=2*\ySep] (x2) {\frame{\input{figures/posets/vertex rooted/vrx2}}};
        \node[xshift=0, yshift=2*\ySep] (xy) {\frame{\input{figures/posets/vertex rooted/vrxy}}};
        \node[xshift=2*\xSep, yshift=2*\ySep] (y2) {\frame{\input{figures/posets/vertex rooted/vry2}}};
        
        \node[xshift=-3*\xSep, yshift=3*\ySep] (x3) {\frame{\input{figures/posets/vertex rooted/vrx3}}};
        \node[xshift=-1*\xSep, yshift=3*\ySep] (x2y) {\frame{\input{figures/posets/vertex rooted/vrx2y}}};
        \node[xshift=1*\xSep, yshift=3*\ySep] (xy2) {\frame{\input{figures/posets/vertex rooted/vrxy2}}};
        \node[xshift=3*\xSep, yshift=3*\ySep] (y3) {\frame{\input{figures/posets/vertex rooted/vry3}}};
        
        \node[xshift=-4*\xSep, yshift=4*\ySep] (x4) {\frame{\input{figures/posets/vertex rooted/vrx4}}};
        \node[xshift=-2*\xSep, yshift=4*\ySep] (x3y) {\frame{\input{figures/posets/vertex rooted/vrx3y}}};
        \node[xshift=0, yshift=4*\ySep] (x2y2) {\frame{\input{figures/posets/vertex rooted/vrx2y2}}};
        \node[xshift=2*\xSep, yshift=4*\ySep] (xy3) {\frame{\input{figures/posets/vertex rooted/vrxy3}}};
        
        \node[xshift=-3*\xSep, yshift=5*\ySep] (x4y) {\frame{\input{figures/posets/vertex rooted/vrx4y}}};
        \node[xshift=1*\xSep, yshift=5*\ySep] (x2y3) {\frame{\input{figures/posets/vertex rooted/vrx2y3}}};
        
        \draw[thick] (N1) -- (x);
        \draw[thick] (N1) -- (y);
        
        \draw[thick] (x) -- (x2);
        \draw[thick] (x) -- (xy);
        \draw[thick] (y) -- (xy);
        \draw[thick] (y) -- (y2);
        
        \draw[thick] (x2) -- (x3);
        \draw[thick] (x2) -- (x2y);
        \draw[thick] (xy) -- (x2y);
        \draw[thick] (xy) -- (xy2);
        \draw[thick] (y2) -- (xy2);
        \draw[thick] (y2) -- (y3);
        
        \draw[thick] (x3) -- (x4);
        \draw[thick] (x3) -- (x3y);
        \draw[thick] (x2y) -- (x3y);
        \draw[thick] (x2y) -- (x2y2);
        \draw[thick] (xy2) -- (x2y2);
        \draw[thick] (xy2) -- (xy3);
        \draw[thick] (y3) -- (xy3);
        
        \draw[thick] (x4) -- (x4y);
        \draw[thick] (x3y) -- (x4y);
        \draw[thick] (x2y2) -- (x2y3);
        \draw[thick] (xy3) -- (x2y3);
    \end{tikzpicture}}

    \captionsetup{width=1.0\linewidth}
  \captionof{figure}{Tri-rooted forests of $\GG$. Marked vertices are denoted by hollow circles and doubly marked vertices are denoted by double circles.}
 \label{fig:runtriiroot}
 \end{center}
\end{figure}


\begin{figure}[H]
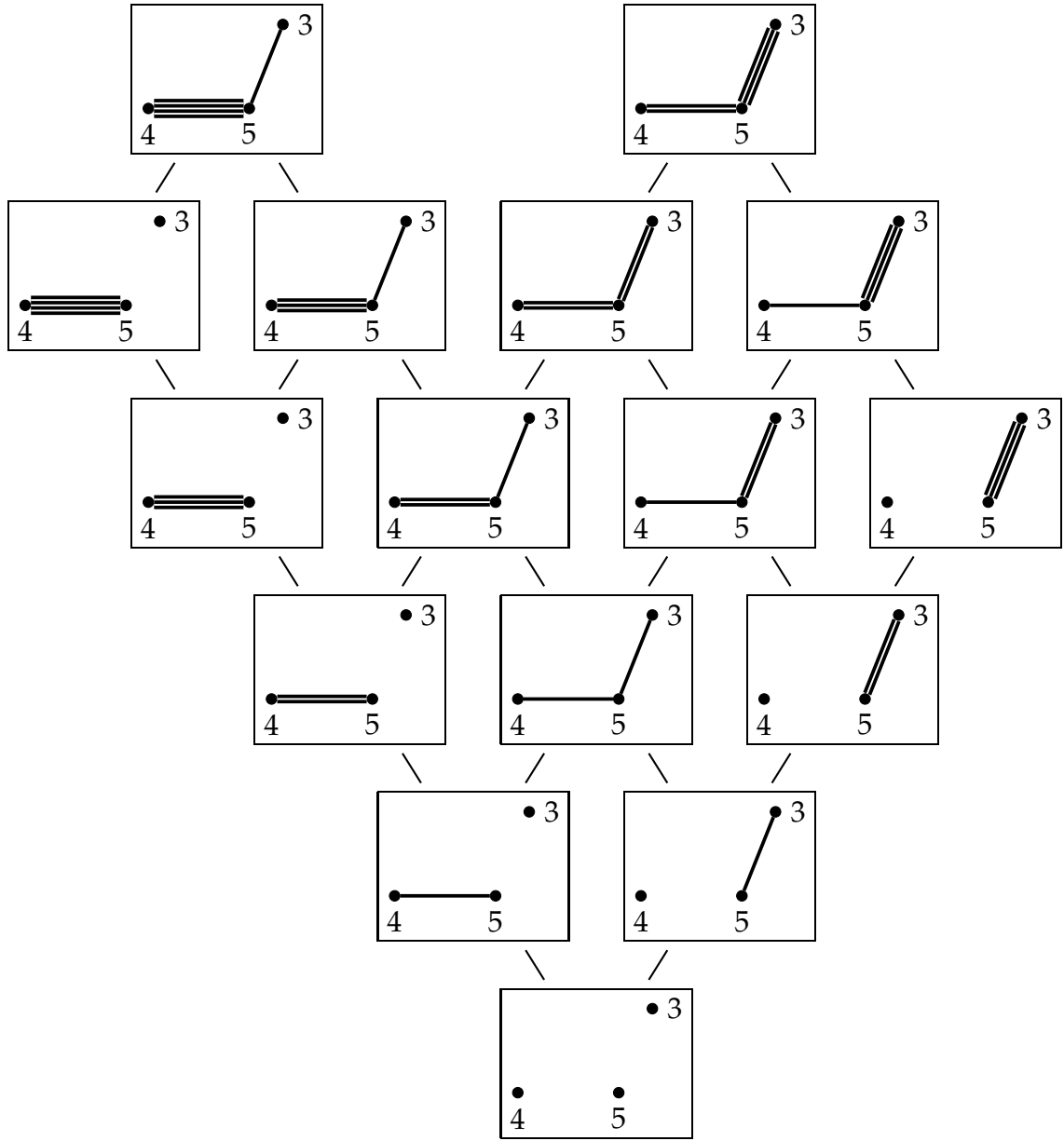

 \begin{center}
 
  \scalebox{1}{        \begin{tikzpicture}[scale=.6]
        
        \node (N1) {\frame{\input{figures/posets/edge rooted/er1}}};
        
        \node[xshift=-1*\xSep, yshift=1*\ySep] (x) {\frame{\input{figures/posets/edge rooted/erx}}};
        \node[xshift=1*\xSep, yshift=1*\ySep] (y) {\frame{\input{figures/posets/edge rooted/ery}}};
        
        \node[xshift=-2*\xSep, yshift=2*\ySep] (x2) {\frame{\input{figures/posets/edge rooted/erx2}}};
        \node[xshift=0, yshift=2*\ySep] (xy) {\frame{\input{figures/posets/edge rooted/erxy}}};
        \node[xshift=2*\xSep, yshift=2*\ySep] (y2) {\frame{\input{figures/posets/edge rooted/ery2}}};
        
        \node[xshift=-3*\xSep, yshift=3*\ySep] (x3) {\frame{\input{figures/posets/edge rooted/erx3}}};
        \node[xshift=-1*\xSep, yshift=3*\ySep] (x2y) {\frame{\input{figures/posets/edge rooted/erx2y}}};
        \node[xshift=1*\xSep, yshift=3*\ySep] (xy2) {\frame{\input{figures/posets/edge rooted/erxy2}}};
        \node[xshift=3*\xSep, yshift=3*\ySep] (y3) {\frame{\input{figures/posets/edge rooted/ery3}}};
        
        \node[xshift=-4*\xSep, yshift=4*\ySep] (x4) {\frame{\input{figures/posets/edge rooted/erx4}}};
        \node[xshift=-2*\xSep, yshift=4*\ySep] (x3y) {\frame{\input{figures/posets/edge rooted/erx3y}}};
        \node[xshift=0, yshift=4*\ySep] (x2y2) {\frame{\input{figures/posets/edge rooted/erx2y2}}};
        \node[xshift=2*\xSep, yshift=4*\ySep] (xy3) {\frame{\input{figures/posets/edge rooted/erxy3}}};
        
        \node[xshift=-3*\xSep, yshift=5*\ySep] (x4y) {\frame{\input{figures/posets/edge rooted/erx4y}}};
        \node[xshift=1*\xSep, yshift=5*\ySep] (x2y3) {\frame{\input{figures/posets/edge rooted/erx2y3}}};
        
        \draw[thick] (N1) -- (x);
        \draw[thick] (N1) -- (y);
        
        \draw[thick] (x) -- (x2);
        \draw[thick] (x) -- (xy);
        \draw[thick] (y) -- (xy);
        \draw[thick] (y) -- (y2);
        
        \draw[thick] (x2) -- (x3);
        \draw[thick] (x2) -- (x2y);
        \draw[thick] (xy) -- (x2y);
        \draw[thick] (xy) -- (xy2);
        \draw[thick] (y2) -- (xy2);
        \draw[thick] (y2) -- (y3);
        
        \draw[thick] (x3) -- (x4);
        \draw[thick] (x3) -- (x3y);
        \draw[thick] (x2y) -- (x3y);
        \draw[thick] (x2y) -- (x2y2);
        \draw[thick] (xy2) -- (x2y2);
        \draw[thick] (xy2) -- (xy3);
        \draw[thick] (y3) -- (xy3);
        
        \draw[thick] (x4) -- (x4y);
        \draw[thick] (x3y) -- (x4y);
        \draw[thick] (x2y2) -- (x2y3);
        \draw[thick] (xy3) -- (x2y3);
    \end{tikzpicture}}
    \captionsetup{width=1.0\linewidth}
  \captionof{figure}{Multicomplex of 3-weighted forests of $\GG$}
 \label{fig:runmulticomp}
 \end{center}
\end{figure}

 
 

\subsection{Further discussions}
As was the case in \cite{biconed}, our approach to proving Stanley's conjecture for the graphic matroid of a triconed graph $G$ centers around choosing an ordering on the edges of $G$, picking a lex minimal tree $B_0$, then using the remaining edges outside of $B_0$ as the variables for the pure multicomplex that realizes the $h$-vector. In \cite{KleeSamper2015} Klee and Samper develop a similar approach to proving Stanley's conjecture (actually a stronger version of it) for all matroids that involves a map from the family of \newword{based matroid} using $B_0^c$ as the variables.  In \cite[Conjecture 3.10]{KleeSamper2015} they suggest the existence of such a map ${\mathcal F}$ that satisfies certain compatibility properties with respect to restriction. Our construction, satisfies the properties $1$ through $4$ they suggested, but we do not know if $5$ is satisfied or not. The argument is pretty much same as what is in Section $6$ of \cite{biconed}. 

We believe next step in extending our methods is for radius $2$ graphs (graphs dominated by a claw) or graphs dominated by a path. When we were only dealing with $3$ special vertices $0,1,2$ the acyclicity of the blueprint could be checked quite easily. But when we have more special vertices, acyclicity of the blueprint gets a lot harder: there can be cycles that involve a number of components. 


Regardless, coming up with an analogue of correctly weighted trees and trirooted forests, and then showing a passivity preserving bijection from the set of original spanning trees does not seem too hard. The really hard part seems to come from trying to come up with an analogue of the map $\phi_2$. The case when all the marks lie on leaves of minimal tree containing all marks within the component seems easy: we simply let each mark empower the adjacent edge within the tree. But when the minimal tree containing all the marks has some marks in the interior (in this paper we only had at most one mark being sandwiched, but there can be more when we have more than $3$ types), it seems pretty hard to decide which edges these marks will empower.

\section*{Acknowledgements}
The research was primarily conducted under $2021$-High school math camp hosted at Texas State University. The authors would like to thank the camp organizers for providing support and a great working environment.

\bibliography{main}
\bibliographystyle{siam}

\end{document}